\documentclass[a4paper,reqno]{amsart}
\usepackage[english]{babel}
\usepackage{amsmath,amssymb,amsthm,exscale,amsopn,fullpage,dsfont,graphicx,mathtools,ragged2e,xspace,yhmath}
\usepackage{varwidth}
\usepackage{xypic}
\usepackage[pdftex,unicode=true,linktocpage,bookmarksopen,hypertexnames=false,
pdfauthor={Lukasz Kubat, Jan Okninski},pdftitle={Irreducible representations of the Chinese monoid}]{hyperref}
\usepackage{memhfixc}

\newtheorem{theorem}{Theorem}[section]
\newtheorem{lemma}[theorem]{Lemma}
\newtheorem{proposition}[theorem]{Proposition}
\newtheorem{corollary}[theorem]{Corollary}

\theoremstyle{definition}
\newtheorem{example}[theorem]{Example}
\newtheorem{remark}[theorem]{Remark}

\DeclareMathOperator{\End}{End}

\DeclareMathOperator{\Ke}{Ker}
\DeclareMathOperator{\Span}{Span}

\newcommand{\diag}[1]{\xymatrix@C=4mm@R=5mm@M=-1mm@W=-1mm{#1}}

\newcommand{\Z}{\mathbb{Z}}
\newcommand{\ol}{\overline}
\newcommand{\arc}{\wideparen}

\author{\L{}ukasz Kubat and Jan Okni\'nski}
\title{Irreducible representations of the Chinese monoid}
\thanks{Research supported by National Science Centre grant 2013/09/B/ST1/04408 (Poland).}

\address{\L{}ukasz Kubat\newline\indent Institute of Mathematics\newline\indent
Warsaw University\newline\indent Banacha 2\newline\indent 02-097 Warsaw, Poland}
\email{\href{mailto:lukasz.kubat@mimuw.edu.pl}{lukasz.kubat@mimuw.edu.pl}}

\address{Jan Okni\'nski\newline\indent Institute of Mathematics\newline\indent
Warsaw University\newline\indent Banacha 2\newline\indent 02-097 Warsaw, Poland}
\email{\href{mailto:okninski@mimuw.edu.pl}{okninski@mimuw.edu.pl}}

\begin{document}

\numberwithin{equation}{section}

\begin{abstract}
All irreducible representations of the Chinese monoid $C_n$, of
any rank $n$, over a nondenumerable algebraically closed field
$K$, are constructed. It turns out that they have a remarkably
simple form and they can be built inductively from irreducible
representations of the monoid $C_2$. The proof shows also that
every such representation is monomial. Since $C_n$ embeds into
the algebra $K[C_n]/J(K[C_n])$, where $J(K[C_n])$
denotes the Jacobson radical of the monoid algebra $K[C_n]$,
a new representation of $C_n$ as a subdirect product of the images
of $C_n$ in the endomorphism algebras of the constructed simple
modules follows.
\end{abstract}

\maketitle

\section{Introduction}

For a positive integer $n$ the monoid $C_n$ defined by the finite
presentation: $C_n =\langle a_1,\dotsc ,a_n\rangle$ with the
defining relations
\begin{equation}\label{relations}
    a_ja_ka_i=a_ka_ja_i=a_ka_ia_j\qquad\text{for }i\le j\le k
\end{equation}
is referred to as the Chinese monoid of rank $n$. It is known that
each element $x$ of $C_n$ has a unique presentation of the form
$x=b_1b_2b_3\dotsm b_n$, where
\begin{align}
    \begin{aligned}
        b_1 & =a_1^{k_{1,1}},\\
        b_2 & =(a_2a_1)^{k_{2,1}} a_2^{k_{2,2}},\\
        b_3 & =(a_3a_1)^{k_{3,1}}(a_3a_2)^{k_{3,2}}a_3^{k_{3,3}},\\
        & \vdotswithin{=}\\
        b_n & =(a_na_1)^{k_{n,1}}(a_na_2)^{k_{n,2}}\dotsm(a_na_{n-1})^{k_{n,n-1}}a_n^{k_{n,n}},
    \end{aligned}\label{canonical}
\end{align}
with all exponents $k_{i,j}$ nonnegative \cite{chinese}. We call
it the canonical form of the element $x\in C_n$. The monoid
algebra $K[C_n]$ over a field $K$, which can be viewed as the
unital algebra defined by the algebra presentation determined by
the relations (\ref{relations}), is called the Chinese algebra of
rank $n$.  The Chinese monoid is related to the so called plactic
monoid, introduced and studied in \cite{las-schut}. Both
constructions are strongly related to Young tableaux, and
therefore to several aspects of representation theory and
algebraic combinatorics. The latter construction became a
classical and powerful tool in representation theory of the full
linear group and in the theory of symmetric functions, via the
Littlewood-Richardson rule (cf.~\cite{fulton}, \cite{las-lec}). It
also plays an important role in quantum groups (in the context of
crystal bases) and in the area of classical Lie
algebras,~\cite{date}, \cite{llt1}, \cite{little}.

The Chinese monoid appeared in the classification of monoids with
the growth function coinciding with that of the plactic monoid
\cite{duchamp}. One of the motivations for a study of the Chinese
monoid is based on an expectation that it might play a
similar role as the plactic monoid in several aspects of
representation theory, quantum algebras, and in algebraic
combinatorics. Combinatorial properties of $C_n$ were studied in
detail in \cite{chinese}. In case $n=2$, the Chinese and the
plactic monoids coincide. The structure of the algebra
$K[C_2]$ is described in \cite{cedo-okn}. In particular, this
algebra is prime and semiprimitive, it is not noetherian and it
does not satisfy any polynomial identity.  For $n=3$ some
information on $K[C_n]$ was obtained in \cite{chin}. In particular
the Jacobson radical of $K[C_3]$ is nonzero, but it is nilpotent,
and the prime spectrum of $K[C_3]$ is pretty well understood. A
surprisingly simple form of the minimal prime ideals of the
algebra $K[C_n]$, for every $n$, was established in \cite{chin2},
\cite{co2}. Namely, every minimal prime ideal $P$ is of the form
$P=\Span_K\{x-y:x,y\in C_n\text{ and }x-y\in P\}$. Hence, in
particular $K[C_n]/P\cong K[C_n/\rho_{P}]$, for the congruence
$\rho_{P}$ on $C_n$ defined by
$\rho_{P}=\{(x,y) \in C_n\times C_n:x-y\in P\}$.
We write $P=I_{\rho_P}$ in this case. It was shown that
every $P$ is generated by a finite set of elements of
the form $x-y$, where $x,y$ are words in the generators
$a_1,\dotsc, a_n$, both of length $2$ or both of length $3$.
Consequently, $K[C_n]/P$ inherits the natural $\Z$-gradation and
this algebra is again defined by a homogeneous semigroup
presentation. In particular, the number of minimal
primes $P$ is finite. Moreover, every $C_n/\rho_{P}$ embeds into
the product $B^k\times\Z^l$, for some
nonnegative integers $k,l$, where $\Z$ is the infinite cyclic group and
$B=\langle p,q \colon qp=1\rangle $ is the bicyclic monoid. The
latter plays an important role in ring theory and in semigroup
theory, \cite{cliff}, \cite{lam}. It was also shown that $C_n$
embeds into the product $\prod_{P}K[C_n]/P$, where $P$ runs over
the set of all minimal primes in $K[C_n]$. Hence $C_n$ embeds into
some $B^r\times \Z^s$. However, the algebra $K[C_n]$ is not
semiprime if $n\geq 3$. Moreover, the description of minimal
primes $P$ of $K[C_n]$ allows to prove that every $K[C_n]/P$ is
semiprimitive and the Jacobson radical of $K[C_n]$ is nilpotent,
and nonzero if $n\ge 3$.

The aim of this paper is to describe all irreducible
representations of $C_n$ over a nondenumerable algebraically
closed field $K$. First, one shows that they are infinite
dimensional unless the dimension is $1$. Then, in our main result,
Theorem~\ref{classification}, all irreducible representations are
constructed. It turns out that they have a remarkably transparent
form. In particular, they can be built inductively from
irreducible representations of the monoid $C_2$, that are
easy to determine. The proof shows that every such representation
is monomial. This is in contrast with the case of
representations of the plactic monoid, as recently shown in
\cite{cedo-kubat-okn}. Since $C_n$ embeds into
$K[C_n]/J(K[C_n])$, where $J(K[C_n])$ denotes the Jacobson
radical of the algebra $K[C_n]$ (see \cite{chin2}), a new
representation of $C_n$ as a subdirect product of the images of
$C_n$ in the endomorphism algebras of the constructed simple
modules also follows.

\section{Background on minimal prime ideals} \label{background}

Throughout the paper, $K$ will stand for a nondenumerable
algebraically closed field, if not stated otherwise. In this
section, we recall from \cite{chin2} the necessary background on
minimal prime ideals of the Chinese algebra $K[C_n]$ of rank
$n\geq 3$.

A finite tree $D$ is associated to $C_n$, whose vertices are diagrams
of certain special type. Each diagram $d$ in $D$ determines a
congruence $\rho (d)$ on the monoid $C_n$ in such a way that the
ideals $I_{\rho(d)}$ corresponding to the leaves $d$ of $D$ are
exactly all the minimal prime ideals of $K[C_n]$.

Each diagram $d$ in $D$ is a graph with $n$ vertices, labeled
$1,\dotsc ,n$ and corresponding to the generators $a_1,\dotsc,a_n$
of $C_n$. For every $d$ in $D$ that is not the root of $D$
there exist $u\le v$, with $u,v\in\{1,\dotsc,n\}$, such that
the vertices $u,\dotsc,v$ are marked (colored black) and the
corresponding generators $a_{u},\dotsc, a_{v}$ are called the used
generators in $d$. Some pairs $k,l$ ($k<l$) of the used generators
can be connected with an edge and then we say that such a pair is
an arc $\arc{a_la_k}$ in $d$. A given generator can be used in at
most one arc. The used (marked) generators not appearing in any
arc are called dots. The root of $D$ is the diagram in which none
of the generators is used. The first level of the tree $D$
(diagrams connected by an edge of $D$ to the root) consists of
$2n-3$ diagrams. There are $n-2$ diagrams with only one of the
generators $a_2,\dotsc, a_{n-1}$ used, and  $n-1$ diagrams with
exactly two consecutive generators $a_s,a_{s+1}$ used in an arc.
Then, if $d$ is in level $t$ of $D$ and it is not a leaf
(by a leaf of the tree $D$ we understand a diagram
containing an arc of the form $\arc{a_ka_1}$ or $\arc{a_na_k}$, for
some $k$) then it is connected by an edge of $D$ to certain diagrams
in level $t+1$ which are obtained from $d$ by adding an arc or
adding a dot, according to the following rules:
\begin{enumerate}
\item If in the last step of construction of the diagram in level $t$ an arc was added,
that is, if the diagram in level $t$ has the form

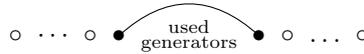
\[\diag{
	\circ & \dotsm & \circ & \bullet \ar@/^1pc/@{-}[rr]
	& {\substack{\text{used}\\ \text{generators}}}  & \bullet & \circ & \dotsc & \circ
}\]
then we can either get, as a diagram in level $t+1$, the diagram\medskip

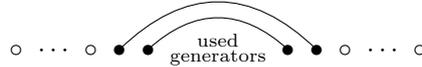
\[\diag{
	\circ & \dotsm & \circ & \bullet \ar@/^1.5pc/@{-}[rrrr] & \bullet \ar@/^1pc/@{-}[rr]
	& {\substack{\text{used}\\ \text{generators}}} & \bullet & \bullet & \circ & \dotsm & \circ
}\]
or one of the following two diagrams

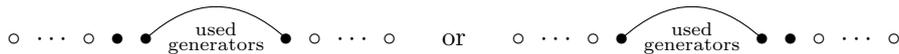
\[\diag{
	\circ & \dotsm & \circ & \bullet & \bullet \ar@/^1pc/@{-}[rr]
	& {\substack{\text{used}\\ \text{generators}}} & \bullet & \circ & \dotsm & \circ}
	\qquad \text{or} \qquad\diag{ \circ & \dotsm & \circ & \bullet \ar@/^1pc/@{-}[rr]
	& {\substack{\text{used}\\ \text{generators}}} & \bullet & \bullet & \circ & \dotsm & \circ
}\]
\item Whereas, if the diagram in level $t$ has the form

\[\diag{
	\circ & \dotsm & \circ & \bullet & {\substack{\text{used}\\ \text{generators}}}
	& \circ & \dotsm & \circ
}\]
then  we can either get, as a diagram in level $t+1$, the diagram

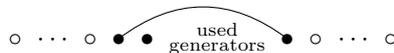
\[\diag{
	\circ & \dotsm & \circ & \bullet \ar@/^1pc/@{-}[rrr] & \bullet
	& {\substack{\text{used}\\ \text{generators}}} & \bullet & \circ & \dotsm & \circ
}\]
or the following diagram

\[\diag{
	\circ & \dotsm & \circ & \bullet & \bullet &
	{\substack{\text{used}\\ \text{generators}}} & \circ & \dotsm & \circ
}\]
\item Similarly, if the diagram in level $t$ has the form

\[\diag{
	\circ & \dotsm & \circ & {\substack{\text{used}\\ \text{generators}}}
	& \bullet & \circ & \dotsm & \circ
}\]
then we can either get, as a diagram in level $t+1$, the diagram

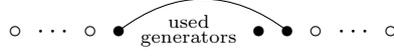
\[\diag{
	\circ & \dotsm & \circ & \bullet \ar@/^1pc/@{-}[rrr] &
	{\substack{\text{used}\\ \text{generators}}} & \bullet & \bullet & \circ & \dotsm & \circ
}\]
or the following diagram

\[\diag{
	\circ & \dotsm & \circ & {\substack{\text{used}\\ \text{generators}}}
	& \bullet & \bullet & \circ & \dotsm & \circ
}\]
\item Finally, after a dot in the first level of $D$
only an arc can be added, so after a diagram

\[\diag{\circ & \dotsm & \circ & \bullet & \circ & \dotsm & \circ } \]
 the following diagram can only occur
 
\[\diag{
	\circ & \dotsm & \circ & \bullet \ar@/^0.4pc/@{-}[rr]
	& \bullet & \bullet & \circ & \dotsm & \circ
}\]
\end{enumerate}
\begin{example}
The diagram in $D$ (for $n=15$) of the form\medskip

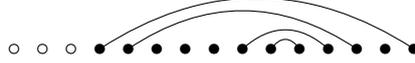
\[\diag{
	\circ & \circ & \circ & \bullet \ar@/^1.6pc/@{-}[rrrrrrrrrrr]
	& \bullet \ar@/^1.2pc/@{-}[rrrrrrrr] & \bullet
	& \bullet & \bullet & \bullet \ar@/^0.6pc/@{-}[rrr]
	& \bullet \ar@/^0.3pc/@{-}[r] & \bullet & \bullet & \bullet & \bullet & \bullet
}\]
arises, in accordance with the rules mentioned above, in the following steps.
First, choose the arc $\arc{a_{11}a_{10}}$ and then the arc $\arc{a_{12}a_9}$.
This leads to the diagram

\[\diag{
	\circ & \circ & \circ & \circ & \circ & \circ & \circ & \circ &
	\bullet \ar@/^0.6pc/@{-}[rrr] & \bullet \ar@/^0.3pc/@{-}[r]
	& \bullet & \bullet & \circ & \circ & \circ
}\]
Next, choose consecutively three dots $a_8,a_7,a_6$.
This yields the following diagram

\[\diag{
	\circ & \circ & \circ & \circ & \circ & \bullet & \bullet & \bullet
	& \bullet \ar@/^0.6pc/@{-}[rrr] & \bullet \ar@/^0.3pc/@{-}[r]
	& \bullet & \bullet & \circ & \circ & \circ
}\]
Then, choosing the arc $\arc{a_{13}a_5}$ we get

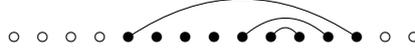
\[\diag{
	\circ & \circ & \circ & \circ & \bullet \ar@/^1.2pc/@{-}[rrrrrrrr]
	& \bullet & \bullet & \bullet & \bullet \ar@/^0.6pc/@{-}[rrr]
	& \bullet \ar@/^0.3pc/@{-}[r] & \bullet & \bullet & \bullet & \circ & \circ
}\]
Finally, choosing the dot $a_{14}$ and then the arc $\arc{a_{15}a_4}$
leads to the considered diagram (which is a leaf).
The full description of $D$ in case $n=3$ and $n=4$ is given in Section~\ref{illustration}.
\end{example}

We shall also consider several homomorphic images of $C_n$ of type
$C_n/\rho$, where $\rho$ is a congruence on $C_n$ generated by certain
pairs of the form $(a_ia_j,a_ja_i)$ and of the form
$(a_ia_ja_k,a_{\sigma(i)}a_{\sigma(j)}a_{\sigma(k)})$ for
some permutations $\sigma$ of $\{i,j,k\}$. Then, for the sake of
simplicity, the image of $a_i$ in $C_n/\rho$ will also be denoted
by $a_i$. Clearly, in monoids of this type we have a degree
function with respect to every generator and we write
$\deg_{a_i}(x)$ for the degree of $x\in C_n/\rho$ in $a_i$.
Moreover, by $\deg(x)$ we mean the total degree of $x$,
that is, $\deg(x)=\sum_{i=1}^n\deg_{a_i}(x)$.

If $u,v\in\{1,\dotsc, n\}$ are such that $u\le v+1$ then we define the monoid
\[C_n^{u,v}=\langle a_1,\dotsc,a_{u-1},a_{v+1},\dotsc ,a_n\rangle\subseteq C_n,\]
which is the Chinese monoid of rank $u-1+n-v$, and its homomorphic image
\[\overline{C}_n^{u,v}=\langle a_1,\dotsc ,a_{u-1},a_{v+1},\dotsc,a_n\rangle
/\big(\begin{smallmatrix}a_1,\dotsc,a_{u-1} & \text{commute}\\
a_{v+1},\dotsc,a_n & \text{commute}\end{smallmatrix}\big),\]
that is, $\ol{C}_n^{u,v}=C_n/\eta$, where $\eta$ is the congruence on $C_n$
generated by all pairs $(a_ia_j,a_ja_i)$ for $i,j<u$ and all pairs $(a_ka_l,a_la_k)$
for $k,l>v$.

By $\Z$ we mean the (multiplicative) infinite cyclic group, with a generator $g$.

If $d_1$ is the diagram with only one used generator $a_s$,
where $1<s<n$, then we associate to it the homomorphism
$\phi_0\colon C_n\longrightarrow\overline{C}_n^{s,s}\times\Z$
defined by
\[\phi_0(a_i) =\begin{cases}
(1,g) & \text{if }i=s,\\
(a_i,1)  & \text{if }i\ne s.
\end{cases}\]
The congruence $\Ke(\phi_0)$ on $C_n$ is
then generated by the pairs:
\begin{alignat*}{2}
	& (a_ia_j,a_ja_i)\qquad && \text{for }i,j\le s,\\
	& (a_ka_l,a_la_k)\qquad &&  \text{for }k,l\ge s.
\end{alignat*}
If $d_1$ is the diagram with only two used generators that form an
arc $\arc{a_{s+1}a_s}$, where $1\le s<n$, then we associate to it the
homomorphism $\psi_0\colon C_n \longrightarrow
\ol{C}_n^{s,s+1}\times B \times\Z$, where
$B=\langle p,q:qp=1\rangle $ is the bicyclic monoid, defined by
\[\psi_0(a_i) = \begin{cases}
    (a_i,p,1)  & \text{if }i<s,\\
    (1,p,g)  & \text{if }i=s,\\
    (1,q,1) &  \text{if }i=s+1,\\
    (a_i,q,1)  & \text{if }i>s+1.
\end{cases}\]
The congruence $\Ke(\psi_0)$ on $C_n$ is then generated by the pairs:
\begin{alignat*}{3}
    & (a_ia_j,a_ja_i),\qquad && (a_ia_{s+1}a_j,a_ja_{s+1}a_i)\qquad && \text{for }i,j\le s,\\
    & (a_ka_l,a_la_k),\qquad && (a_ka_sa_l,a_la_sa_k)\qquad && \text{for }k,l>s.
\end{alignat*}
Now, we define \[\kappa(d_1)\colon C_n\longrightarrow\ol{C}_n^{u,v}\times S_1\]
as $\kappa(d_1)=\phi_0$ (and then $(u,v)=(s,s)$ and $S_1=\Z$) in case $d_1$
is the diagram with only one used generator $a_s$, or $\kappa(d_1)=\psi_0$
(and then $(u,v)=(s,s+1)$ and $S_1=B\times\Z$) in case $d_1$ is the diagram
with only two used generators that form an arc $\arc{a_{s+1}a_s}$.
Moreover, let $\rho(d_1)=\Ke\kappa(d_1)$.

So, the homomorphisms and congruences described above are associated
to the $2n-3$ diagrams from the first level of $D$. The procedure described
below allows us to associate (inductively) a homomorphism
$\kappa(d)\colon C_n\longrightarrow\ol{C}_n^{u,v}\times(B\times\Z)^k\times\Z^l$,
where $u,v$ and $k,l$ depend on $d$, and the congruence $\rho(d)=\Ke\kappa(d)$
to every diagram $d$ at the level $>1$ of $D$. However, in contrast to the congruences
from the first level, the generators of $\rho(d)$ are much
harder to determine, see \cite{co2}.

Assume that a diagram $d_t$ in level $t\ge 1$ of the tree $D$ has
been constructed and it is not a leaf. Assume also that the homomorphism
$\kappa(d_t)\colon C_n\longrightarrow\ol{C}_n^{u,v}\times S_t$, where
$S_t=(B\times\Z)^k\times\Z^l$, together with the congruence
$\rho(d_t)=\Ke\kappa(d_t)$ have been defined. Here, the indices $u,v$ correspond
to the used generators $a_u,\dotsc,a_v$ in the diagram $d_t$, whereas
the nonnegative integers $k,l$ correspond to the number of arcs and dots,
respectively, used in the construction of the diagram $d_t$.
Moreover, let $d_{t+1}$ be a diagram at the level $t+1$
of $D$ that is connected to $d_t$ by an edge in $D$.

If $d_{t+1}$
is obtained by adding a dot to $d_t$ (then this is either $a_{u-1}$
or $a_{v+1}$), we have a homomorphism
\[\phi_t\colon\ol{C}_n^{u,v}\times S_t\longrightarrow\ol{C}_n^{u',v'} \times\Z\times S_t\]
given by
\[\phi_t(a_i,x)=\begin{cases}
    (1,g,x) & \text{if }i=s,\\
    (a_i,1,x) & \text{if }i<u\text{ or }i>v\text{ but }i\ne s,
\end{cases}\]
where $s=u-1$ (and then $(u',v')=(u-1,v)$) or $s=v+1$ (and then
$(u',v')=(u,v+1)$), depending on which of the two possible dots was
added.

If $d_{t+1}$ is obtained by adding an arc to $d_t$ (then this arc is
$\arc{a_{v+1}a_{u-1}}$), we have a homomorphism
\[\psi_t \colon\ol{C}_n^{u,v}\times S_t\longrightarrow\ol{C}_n^{u-1, v+1} \times B \times\Z\times S_t\]
given by
\[\psi_t(a_i,x) =\begin{cases}
    (a_i,p,1,x) & \text{if }i<u-1,\\
    (1,p,g,x) & \text{if }i=u-1,\\
    (1,q,1,x)  & \text{if }i=v+1,\\
    (a_i,q,1,x) & \text{if }i>v+1,
\end{cases}\]
and then we put $(u',v')=(u-1,v+1)$. Furthermore, we put $S_{t+1}=\Z\times S_t$
in case $d_{t+1}$ is obtained by adding a dot to $d_t$, and $S_{t+1}=B\times\Z\times S_t$
in case $d_{t+1}$ is obtained by adding an arc to $d_t$.

Now, we define \[\kappa(d_{t+1})\colon C_n\longrightarrow\ol{C}_n^{u',v'}\times S_{t+1}\]
as a composition
\[\kappa(d_{t+1})=\begin{cases}
	\phi_t\circ\kappa(d_t) & \text{ if }d_{t+1}\text{ is obtained from }d_t\text{ by adding a dot},\\
	\psi_t\circ\kappa(d_t) & \text{ if }d_{t+1}\text{ is obtained from }d_t\text{ by adding an arc}.
\end{cases}\]
Moreover, let $\rho(d_{t+1})=\Ke\kappa(d_{t+1})$. Then, of course, $\rho(d_t)\subseteq\rho(d_{t+1})$.

Summarizing, if $d_0,d_1,\dotsc,d_m=d$ in a branch in $D$ then
$\rho(d_0)\subseteq\rho(d_1)\subseteq\dotsb\subseteq\rho(d_m)=\rho(d)$
(by $\rho(d_0)$, for the root $d_0$ of $D$, we mean the trivial congruence on $C_n$).
However, if $m>1$, then as was mentioned at the beginning of this section,
generators of the congruence $\rho(d)=\Ke\kappa(d)$ are hard to determine explicitly.
Though, if the diagram $d$ is of some special shape (e.g. one of the shapes listed below),
then looking at the embedding
$C_n/\rho(d)\longrightarrow\ol{C}_n^{u,v}\times(B\times\Z)^k\times\Z^l$
(for some $u,v$ and some $k,l$), induced by the homomorphism $\kappa(d)$,
it is quite easy to derive some relations that must hold in $C_n/ \rho(d)$ and which will be needed later.

Namely, if $d$ is a diagram of the form

\[ \diag{ \circ & \dotsm & \circ & \underset{s}{\bullet} & \circ & \dotsm & \circ }\]
consisting of a single dot $a_s$, then the following equalities hold in $C_n/\rho(d)$:
\begin{alignat}{2}
	& a_ia_j = a_ja_i\qquad & & \text{for all }i,j\le s, \label{rel1-1}\\
	& a_ka_l=a_la_k\qquad & & \text{for all }k,l\ge s. \label{rel1-2}
\end{alignat}
If $d$ is a diagram of the form

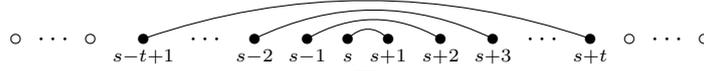
\[\diag{
	\circ & \dotsm & \circ & \underset{s-t+1}{\bullet} \ar@/^1.2pc/@{-}[rrrrrrrrr] &
	\dotsm & \underset{s-2}{\bullet} \ar@/^0.9pc/@{-}[rrrrr] & \underset{s-1}{\bullet} \ar@/^0.6pc/@{-}[rrr]
        & \underset{s}{\bullet} \ar@/^0.3pc/@{-}[r] & \underset{s+1}{\bullet} & \underset{s+2}{\bullet}
        & \underset{s+3}{\bullet} & \dotsm & \underset{s+t}{\bullet} & \circ & \dotsm & \circ
}\]
consisting of $t>0$ consecutive arcs $\arc{a_{s+1}a_s},\dotsc,\arc{a_{s+t}a_{s-t+1}}$,
then the following equalities hold in $C_n/\rho(d)$:
\begin{alignat}{2}
	& a_ia_j = a_ja_i\qquad & & \text{for all }i,j\le s,\label{rel2-1} \\
	& a_ka_l=a_la_k\qquad & & \text{for all }k,l>s,\label{rel2-2} \\
	& a_ia_{s+r}a_j=a_ja_{s+r}a_i\qquad & & \text{for all }i,j\le s-r+1\text{, where }r=1,\dotsc,t,\label{rel2-3} \\
	& a_ka_{s-r+1}a_l=a_la_{s-r+1}a_k\qquad & & \text{for all }k,l\ge s+r\text{, where }r=1,\dotsc,t.\label{rel2-4}
\end{alignat}
If $d$ is a diagram of the form

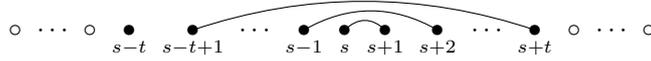
\[\diag{
	\circ & \dotsm & \circ & \underset{s-t}{\bullet} & \underset{s-t+1}{\bullet} \ar@/^0.9pc/@{-}[rrrrrrr]
	& \dotsm & \underset{s-1}{\bullet} \ar@/^0.6pc/@{-}[rrr] & \underset{s}{\bullet} \ar@/^0.3pc/@{-}[r]
	& \underset{s+1}{\bullet} & \underset{s+2}{\bullet} & \dotsm & \underset{s+t}{\bullet} & \circ & \dotsm & \circ
}\]
consisting of $t>0$ consecutive arcs $\arc{a_{s+1}a_s},\dotsc,\arc{a_{s+t}a_{s-t+1}}$
and a single dot $a_{s-t}$, then the following equalities hold in $C_n/\rho(d)$:
\begin{alignat}{2}
	& a_ia_j = a_ja_i\qquad & & \text{for all }i,j\le s, \label{rel3-1} \\
	& a_ka_l=a_la_k\qquad & & \text{for all }k,l>s, \label{rel3-2} \\
	& a_ka_{s-t}a_l=a_la_{s-t}a_k & & \text{for all }k,l\ge s+t, \label{rel3-3} \\
	& a_ia_{s+r}a_j=a_ja_{s+r}a_i\qquad & & \text{for all }i,j\le s-r+1\text{, where }r=1,\dotsc,t,\label{rel3-4} \\
	& a_ka_{s-r+1}a_l=a_la_{s-r+1}a_k\qquad & & \text{for all }k,l\ge s+r\text{, where }r=1,\dotsc,t.\label{rel3-5}
\end{alignat}
Similarly, if $d$ is a diagram of the form

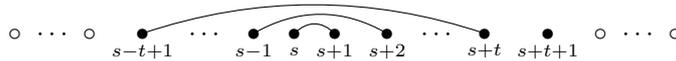
\[\diag{
	\circ & \dotsm & \circ & \underset{s-t+1}{\bullet} \ar@/^0.9pc/@{-}[rrrrrrr]
	& \dotsm & \underset{s-1}{\bullet} \ar@/^0.6pc/@{-}[rrr] & \underset{s}{\bullet} \ar@/^0.3pc/@{-}[r]
	& \underset{s+1}\bullet & \underset{s+2}{\bullet} & \dotsm & \underset{s+t}{\bullet}
	& \underset{s+t+1}{\bullet} & \circ & \dotsm & \circ
}\]
consisting of $t>0$ consecutive arcs $\arc{a_{s+1}a_s},\dotsc,\arc{a_{s+t}a_{s-t+1}}$
and a single dot $a_{s+t+1}$, then equalities dual to (\ref{rel3-1})--(\ref{rel3-5}) also
can be derived. However, these equalities will not be used explicitly in the paper.

\section{Irreducible representations}\label{irr}

Our first result shows that infinite dimensional
simple $K[C_n]$-modules will be crucial.

\begin{proposition} \label{findim}
Let $\phi\colon C_n\longrightarrow\End_K(V)$ be an irreducible
representation of $C_n$ over a field $K$.
\begin{enumerate}
	\item If $\phi (a_na_1)=0$ then either $\phi (a_n)=0$ or $\phi(a_1)=0$.
	\item If $\dim_KV<\infty$, then $\phi (C_n)$ is commutative,
	hence $\dim_KV=1$ if $K$ is algebraically closed.
\end{enumerate}
\begin{proof}
First, we claim that $a_nC_na_1 \subseteq a_na_1C_n$.
Let $x\in C_n$. From the canonical form of elements in $C_n$
it follows that $xa_1=x_j\dotsm x_n$ for some $j\in\{1,\dotsc,n\}$,
where $x_j, \dotsc , x_n \in C_n$ and $x_j\in a_1C_n$ if $j=1$,
and $x_j\in a_ja_1C_n$ if $j>1$.
Since $a_na_ja_1=a_na_1a_j$, the claim follows.

Suppose that $\phi (a_na_1)=0$. Then
$\phi (a_n)\phi(C_n)\phi(a_1)=\phi(a_n)\phi(a_1)\phi(C_n)=0$.
Since $\phi$ is irreducible, it follows that either
$\phi (a_n)=0 $ or $\phi(a_1)=0$.

In order to prove the second assertion, consider $z=a_na_1\in C_n$.
If $\phi(z)=0$ then, by the first part of the proof, $\phi$ comes from
an irreducible representation of $C_{n-1}$. Hence, the result
follows by induction in this case. Otherwise, $\phi(z)\ne0$ is an
invertible element in the simple algebra $R=\Span_K\phi(C_n)$,
because it is central. Thus, in particular, $\phi(a_n)$ is invertible
in $R$, so the relations (\ref{relations}) defining $C_n$ easily imply that
$\phi(C_n)$ is commutative and the assertion follows.
\end{proof}
\end{proposition}

Our next aim is to construct a family of simple left
$K[C_n]$-modules in case $n$ is even. Later we shall see that
these modules are of special interest, because they are the corner
stone of an inductive classification of all simple left modules
over the algebra $K[C_n]$.
\begin{proposition}\label{simple}
Let $V$ be a $K$-linear space with basis
$\{e_{i_1,\dotsc,i_s}:i_1,\dotsc,i_s\ge 0\}$ for some $s\ge 1$.
Moreover, let $0\ne\lambda_1,\dotsc,\lambda_s\in K$ and $n=2s$.
Then the action of $a_1,\dotsc,a_n\in C_n$ on $V$ defined by
\[a_je_{i_1,\dotsc,i_s}=\begin{cases}
	\lambda_je_{i_1,\dotsc,i_{j-1},i_j+1,\dotsc,i_s+1} & \text{if }j\le s,\\
	e_{i_1,\dotsc,i_{n-j},i_{n-j+1}-1,\dotsc,i_s-1} & \text{if }j>s\text{ and }i_k>0\text{ for all }k>n-j,\\
	0 & \text{if }j>s\text{ and }i_k=0\text{ for some }k>n-j
\end{cases}\]
makes $V=V(\lambda_1,\dotsc,\lambda_s)$ a simple left $K[C_n]$-module.
Moreover, if $0\ne\mu_1,\dotsc,\mu_s\in K$ then we have
$V(\lambda_1,\dotsc,\lambda_s)\cong V(\mu_1,\dotsc,\mu_s)$ as left
$K[C_n]$-modules if and only if $\lambda_i=\mu_i$ for all $i=1,\dotsc,s$.
\begin{proof}
First, we have to check that the defined action of $a_1,\dotsc,a_n\in C_n$
on $V$ respects the Chinese relations. So, we have to prove that
\begin{align*}
	(a_la_ka_j-a_la_ja_k)V & =a_l(a_ka_j-a_ja_k)V=0
	\intertext{and}
	(a_la_ka_j-a_ka_la_j)V & =(a_la_k-a_ka_l)a_jV=0
\end{align*}
for all $j\le k\le l$. Since we have $(a_ka_j-a_ja_k)V=0$
for all $j,k\le s$ and $(a_la_k-a_ka_l)V=0$ for all $k,l>s$, it is enough to show that:
\begin{enumerate}
	\item $(a_la_ka_j-a_la_ja_k)V=0$ for all $j\le s<k\le l$ such that $j+l\le n$,
	\item $(a_la_ka_j-a_la_ja_k)V=0$ for all $j\le s<k\le l$ such that $j+k\le n<j+l$,
	\item $(a_la_ka_j-a_la_ja_k)V=0$ for all $j\le s<k\le l$ such that $n<j+k$,
	\item $(a_la_ka_j-a_ka_la_j)V=0$ for all $j\le k\le s<l$ such that $k+l\le n$,
	\item $(a_la_ka_j-a_ka_la_j)V=0$ for all $j\le k\le s<l$ such that $j+l\le n<k+l$,
	\item $(a_la_ka_j-a_ka_la_j)V=0$ for all $j\le k\le s<l$ such that $j+l>n$.
\end{enumerate}
It is easy to verify that we have, respectively:
\begin{enumerate}
	\item If $k<l$ then
	\[a_la_ka_je_{i_1,\dotsc,i_s}=a_la_ja_ke_{i_1,\dotsc,i_s}=\begin{cases}
	\lambda_j e_{i_1,\dotsc,i_{j-1},i_j+1,\dotsc,i_{n-l}+1,i_{n-l+1},\dotsc,i_{n-k},i_{n-k+1}-1,\dotsc,i_s-1}\\
	\hfill\text{if }i_p>0\text{ for all }p>n-k,\\
	0\hfill\text{otherwise}.\end{cases}\]
	Whereas, if $k=l$ then
	\[a_la_ka_je_{i_1,\dotsc,i_s}=a_la_ja_ke_{i_1,\dotsc,i_s}=\begin{cases}
	\lambda_j e_{i_1,\dotsc,i_{j-1},i_j+1,\dotsc,i_{n-k}+1,i_{n-k+1}-1,\dotsc,i_s-1}\\
	\hfill\text{if }i_p>0\text{ for all }p>n-k,\\
	0\hfill\text{otherwise}.\end{cases}\]
	\item If $j+l>n+1$ then
	\[a_la_ka_je_{i_1,\dotsc,i_s}=a_la_ja_ke_{i_1,\dotsc,i_s}=\begin{cases}
	\lambda_j e_{i_1,\dotsc,i_{n-l},i_{n-l+1}-1,\dotsc,i_{j-1}-1,i_j,\dotsc,i_{n-k},i_{n-k+1}-1,\dotsc,i_s-1}\\
	\hfill\text{if }i_p>0\text{ for all }n-l<p<j,\\
	\hfill\text{and }i_q>0\text{ for all }q>n-k,\\
	0\hfill\text{otherwise}.\end{cases}\]
	Whereas, if $j+l=n+1$ then
	\[a_la_ka_je_{i_1,\dotsc,i_s}=a_la_ja_ke_{i_1,\dotsc,i_s}=\begin{cases}
	\lambda_j e_{i_1,\dotsc,i_{n-k},i_{n-k+1}-1,\dotsc,i_s-1}\text{ if }i_p>0\text{ for all }p>n-k,\\
	0\hfill\text{otherwise}.\end{cases}\]
	\item If $k<l$ and $j+k>n+1$ then
	\[a_la_ka_je_{i_1,\dotsc,i_s}=a_la_ja_ke_{i_1,\dotsc,i_s}=\begin{cases}
	\lambda_j e_{i_1,\dotsc,i_{n-l},i_{n-l+1}-1,\dotsc,i_{n-k}-1,i_{n-k+1}-2,\dotsc,i_{j-1}-2,i_j-1,\dotsc,i_s-1}\\
	\hfill\text{if }i_p>0\text{ for all }n-l<p\le n-k,\\
	\hfill i_q>1\text{ for all }n-k<q<j,\\
	\hfill\text{and }i_r>0\text{ for all }r\ge j,\\
	0\hfill\text{otherwise}.\end{cases}\]
	If $k=l$ and $j+k>n+1$ then
	\[a_la_ka_je_{i_1,\dotsc,i_s}=a_la_ja_ke_{i_1,\dotsc,i_s}=\begin{cases}
	\lambda_j e_{i_1,\dotsc,i_{n-k},i_{n-k+1}-2,\dotsc,i_{j-1}-2,i_j-1,\dotsc,i_s-1}\\
	\hfill\text{if }i_p>1\text{ for all }n-k<p<j,\\
	\hfill\text{and }i_q>0\text{ for all }q\ge j,\\
	0\hfill\text{otherwise}.\end{cases}\]
	Whereas, if $j+k=n+1$ then
	\[a_la_ka_je_{i_1,\dotsc,i_s}=a_la_ja_ke_{i_1,\dotsc,i_s}=\begin{cases}
	\lambda_j e_{i_1,\dotsc,i_{n-l},i_{n-l+1}-1,\dotsc,i_s-1}\text{ if }i_p>0\text{ for all }p>n-l,\\
	0\hfill\text{otherwise}.\end{cases}\]
	\item If $j<k$ then
	\[a_la_ka_je_{i_1,\dotsc,i_s}=a_ka_la_je_{i_1,\dotsc,i_s}=
	\lambda_j\lambda_ke_{i_1,\dotsc,i_{j-1},i_j+1,\dotsc,i_{k-1}+1,i_k+2,\dotsc,i_{n-l}+2,i_{n-l+1}+1,\dotsc,i_s+1}.\]
	Whereas, if $j=k$ then
	\[a_la_ka_je_{i_1,\dotsc,i_s}=a_ka_la_je_{i_1,\dotsc,i_s}=\lambda_j^2e_{i_1,\dotsc,i_{j-1},i_j+2,\dotsc,i_{n-l}+2,i_{n-l+1}+1,\dotsc,i_s+1}.\]
	\item If $k+l>n+1$ then
	\[a_la_ka_je_{i_1,\dotsc,i_s}=a_ka_la_je_{i_1,\dotsc,i_s}=
	\lambda_j\lambda_ke_{i_1,\dotsc,i_{j-1},i_j+1,\dotsc,i_{n-l}+1,i_{n-l+1},\dotsc,i_{k-1},i_k+1,\dotsc,i_s+1}.\]
	Whereas, if $k+l=n+1$ then
	\[a_la_ka_je_{i_1,\dotsc,i_s}=a_ka_la_je_{i_1,\dotsc,i_s}=
	\lambda_j\lambda_ke_{i_1,\dotsc,i_{j-1},i_j+1,\dotsc,i_s+1}.\]
	\item If $j<k$ and $j+l>n+1$ then
	\[a_la_ka_je_{i_1,\dotsc,i_s}=a_ka_la_je_{i_1,\dotsc,i_s}=\begin{cases}
	\lambda_j\lambda_k e_{i_1,\dotsc,i_{n-l},i_{n-l+1}-1,\dotsc,i_{j-1}-1,i_j,\dotsc,i_{k-1},i_k+1,\dotsc,i_s+1}\\
	\hfill\text{if }i_p>0\text{ for all }n-l<p<j,\\
	0\hfill\text{otherwise}.\end{cases}\]
	If $j=k$ and $j+l>n+1$ then
	\[a_la_ka_je_{i_1,\dotsc,i_s}=a_ka_la_je_{i_1,\dotsc,i_s}=\begin{cases}
	\lambda_j^2 e_{i_1,\dotsc,i_{n-l},i_{n-l+1}-1,\dotsc,i_{j-1}-1,i_j+1,\dotsc,i_s+1}\\
	\hfill\text{if }i_p>0\text{ for all }n-l<p<j,\\
	0\hfill\text{otherwise}.\end{cases}\]
	Whereas, if $j+l=n+1$ then
	\[a_la_ka_je_{i_1,\dotsc,i_s}=a_ka_la_je_{i_1,\dotsc,i_s}=
	\lambda_j\lambda_k e_{i_1,\dotsc,i_{k-1},i_k+1,\dotsc,i_s+1}.\]
\end{enumerate}

Now, let us prove that the $K[C_n]$-module $V$ is simple. First,
it can be easily verified that for each $j<s$ we have
$a_{n-j}a_je_{i_1,\dotsc,i_s}=\lambda_je_{i_1,\dotsc,i_{j-1},i_j+1,i_{j+1},\dotsc,i_s}$
(that is, the action of $a_{n-j}a_j$ on $e_{i_1,\dotsc,i_s}$
increases the index $i_j$ by one and leaves other indices
unchanged) and
$a_se_{i_1,\dotsc,i_s}=\lambda_se_{i_1,\dotsc,i_{s-1},i_s+1}$
(that is, the action of $a_s$ on $e_{i_1,\dotsc,i_s}$ increases
the index $i_s$ by one and leaves other indices unchanged).
Therefore,
\[e_{i_1,\dotsc,i_s}=(\lambda_1^{-1}a_{n-1}a_1)^{i_1}\dotsm
(\lambda_{s-1}^{-1}a_{s+1}a_{s-1})^{i_{s-1}}(\lambda_s^{-1}a_s)^{i_s}e_{0,\dotsc,0}\]
for all $i_1,\dotsc,i_s\ge 0$. Hence, to prove simplicity of $V$, it suffices to check that
$e_{0,\dotsc,0}\in K[C_n]v$ for each $0\ne v\in V$. So, let
$0\ne v=\sum_{i_1,\dotsc,i_s=0}^r\lambda_{i_1,\dotsc,i_s}e_{i_1,\dotsc,i_s}$ for some
$\lambda_{i_1,\dotsc,i_s}\in K$ be fixed. Then define
\begin{align*}
	m_1 & =\max\{i_1:\lambda_{i_1,\dotsc,i_s}\ne 0\text{ for some }i_2,\dotsc,i_s\},\\
	m_2 & =\max\{i_2:\lambda_{m_1,i_2,\dotsc,i_s}\ne 0\text{ for some }i_3,\dotsc,i_s\},\\
	m_3 & =\max\{i_3:\lambda_{m_1,m_2,i_3,\dotsc,i_s}\ne 0\text{ for some }i_4,\dotsc,i_s\},\\
	& \vdotswithin{=}\\
	m_s & =\max\{i_s:\lambda_{m_1,\dotsc,m_{s-1},i_s}\ne 0\}.
\end{align*}
Because, for each $j<s$, we have
\[a_{n-j+1}a_{j+1}e_{i_1,\dotsc,i_s}=\begin{cases}
\lambda_{j+1}e_{i_1,\dotsc,i_{j-1},i_j-1,i_{j+1},\dotsc,i_s} & \text{if }i_j>0,\\
0 & \text{if }i_j=0\end{cases}\]
(that is, the action of $a_{n-j+1}a_{j+1}$ on $e_{i_1,\dotsc,i_s}$
decreases the index $i_j$ by one, if possible, and leaves other
indices unchanged) and because
\[a_{s+1}e_{i_1,\dotsc,i_s}=\begin{cases}e_{i_1,\dotsc,i_{s-1},i_s-1}
& \text{if }i_s>0,\\ 0 & \text{if }i_s=0\end{cases}\]
(that is, the action of $a_{s+1}$ on $e_{i_1,\dotsc,i_s}$
decreases the index $i_s$ by one, if possible, and leaves other
indices unchanged), we conclude that
\[e_{0,\dotsc,0}=(\lambda_{m_1,\dotsc,m_s}\lambda_2^{m_1}\dotsm\lambda_s^{m_{s-1}})^{-1}
a_{s+1}^{m_s}(a_{s+2}a_s)^{m_{s-1}}\dotsm(a_na_2)^{m_1}v\in K[C_n]v,\]
as claimed.

Finally, note that isomorphic modules have equal annihilators
and $(a_{n-i+1}a_i-\lambda_i)V=0$ for each $i=1,\dotsc,s$.
Hence, if also $(a_{n-i+1}a_i-\mu_i)V=0$ for some $i\in\{1,\dotsc,s\}$,
then $(\lambda_i-\mu_i)V=0$ and, in consequence,
$\lambda_i=\mu_i$. Thus the last part of the proposition also follows.
\end{proof}
\end{proposition}

It is worth to note that the modules constructed in Proposition~\ref{simple}
can be obtained by a successive application of the construction presented
in Proposition~\ref{inductive}, starting with the left $K[C_2]$-module $Z$
with basis $\{e_i:i\ge 0\}$, and with the action of $a_1,a_2\in C_2$
on $Z$ defined by
\[a_1e_i=\lambda_1 e_{i+1},\qquad a_2e_i=\begin{cases}e_{i-1}
& \text{if }i>0,\\ 0 & \text{if }i=0\end{cases}\]
for some $0\ne\lambda_1\in K$.
(Notice that such a module $Z$ is a straightforward generalization of
the classical simple $K[B]$-module, considered for example in \cite[p.~195 and Ex.~11.9.]{lam}.)
This is fully explained below.

\begin{proposition}\label{inductive}
Let $U$ be a left $K[C_n^{s,s+1}]$-module with basis
$\{e_{i_1,\dotsc,i_{s-1}}:i_1,\dotsc,i_{s-1}\ge 0\}$, where $n=2s$ for some $s\ge 1$.
Assume that $(a_ia_j-a_ja_i)U=(a_ka_l-a_la_k)U=0$ for all $i,j<s$ and
$k,l>s+1$. Assume also that for each $0\ne u\in U$ and
$i_1,\dotsc,i_{s-1}\ge 0$ there exists $x\in K\cdot C_n^{s,s+1}$
satisfying $xu=e_{i_1,\dotsc,i_{s-1}}$. Moreover, let $V$ be a
$K$-linear space with basis $\{f_i:i\ge 0\}$. Then, for each
$0\ne\lambda_s\in K$, the action of $a_1,\dotsc,a_n\in C_n$ on the
$K$-linear space $W=U\otimes_KV$ with basis
$\{e_{i_1,\dotsc,i_s}=e_{i_1,\dotsc,i_{s-1}}\otimes f_{i_s}:i_1,\dotsc,i_s\ge 0\}$
defined by
\[a_je_{i_1,\dotsc,i_s}=\begin{cases}
a_je_{i_1,\dotsc,i_{s-1}}\otimes f_{i_s+1} & \text{if }j<s,\\
\lambda_s e_{i_1,\dotsc,i_{s-1}}\otimes f_{i_s+1} & \text{if }j=s,\\
e_{i_1,\dotsc,i_{s-1}}\otimes f_{i_s-1} & \text{if }j=s+1\text{ and }i_s>0,\\
a_je_{i_1,\dotsc,i_{s-1}}\otimes f_{i_s-1} & \text{if }j>s+1\text{ and }i_s>0,\\
0 & \text{otherwise}\end{cases}\]
makes $W$ a left $K[C_n]$-module such that
$(a_ia_j-a_ja_i)W=(a_ka_l-a_la_k)W=0$ for all $i,j\le s$ and
$k,l>s$. Moreover, for each $0\ne w\in W$ and $i_1,\dotsc,i_s\ge 0$
there exists $x\in K\cdot C_n$ satisfying $xw=e_{i_1,\dotsc,i_s}$.
In particular, $W$ is a simple left $K[C_n]$-module.
\begin{proof}
First, it is convenient to define
$a_se_{i_1,\dotsc,i_{s-1}}=\lambda_s e_{i_1,\dotsc,i_{s-1}}$ and
$a_{s+1}e_{i_1,\dotsc,i_{s-1}}=e_{i_1,\dotsc,i_{s-1}}$ for all
$i_1,\dotsc,i_{s-1}\ge 0$ . With this notation the action of
$a_1,\dotsc,a_n\in C_n$ on the basis of $W$ can be rewritten as
\[a_je_{i_1,\dotsc,i_s}=\begin{cases}a_je_{i_1,\dotsc,i_{s-1}}\otimes f_{i_s+1} & \text{if }j\le s,\\
a_je_{i_1,\dotsc,i_{s-1}}\otimes f_{i_s-1} & \text{if }j>s\text{ and }i_s>0,\\
0 & \text{if }j>s\text{ and }i_s=0.\end{cases}\]
It is almost obvious that the defined action respects the Chinese
relations not involving $a_s$ and $a_{s+1}$. Moreover, if $i,j\le s$ then
\[(a_ia_j-a_ja_i)e_{i_1,\dotsc,i_s}=(a_ia_j-a_ja_i)e_{i_1,\dotsc,i_{s-1}}\otimes f_{i_s+2}=0,\]
hence $(a_ia_j-a_ja_i)W=0$. Similarly, if $k,l>s$ then
\[(a_ka_l-a_la_k)e_{i_1,\dotsc,i_s}=\begin{cases}(a_ka_l-a_la_k)e_{i_1,\dotsc,i_{s-1}}
\otimes f_{i_s-2} & \text{if }i_s>1,\\ 0 & \text{if }i_s\le 1\end{cases}=0,\]
hence $(a_ka_l-a_la_k)W=0$. Since $(a_{s+1}a_s-\lambda_s)W=0$,
to prove that $W$ is indeed a $K[C_n]$-module it suffices to check that:
\begin{enumerate}
	\item $(a_ja_{s+1}a_i-a_{s+1}a_ja_i)W=0$ for all $i\le j\le s$,
	\item $(a_ja_{s+1}a_i-a_ja_ia_{s+1})W=0$ for all $i\le s<j$,
	\item $(a_sa_ja_i-a_ja_sa_i)W=0$ for all $i\le s<j$,
	\item $(a_ja_ia_s-a_ja_sa_i)W=0$ for all $s<i\le j$,
	\item $(a_{s+1}a_{s+1}a_i-a_{s+1}a_ia_{s+1})W=0$ for all $i\le s$,
	\item $(a_sa_ia_s-a_ia_sa_s)W=0$ for all $i>s$.
\end{enumerate}
Let $W_0$ denote the subspace of $W$ spanned by the set
$\{e_{i_1,\dotsc,i_s}:i_1,\dotsc,i_{s-1}\ge 0\text{ and }i_s=0\}$,
and $W_+$ denote the subspace of $W$ spanned by
the set $\{e_{i_1,\dotsc,i_s}:i_1,\dotsc,i_{s-1}\ge 0\text{ and }i_s>0\}$.
Then we have, respectively:
\begin{enumerate}
	\item $(a_ja_{s+1}a_i-a_{s+1}a_ja_i)W=(a_ja_{s+1}-a_{s+1}a_j)a_iW=0$,
	because $a_iW\subseteq W_+$ for $i\le s$ and $(a_ja_{s+1}-a_{s+1}a_j)W_+=0$ for $j\le s$.
	\item $(a_ja_{s+1}a_i-a_ja_ia_{s+1})W=a_j(a_{s+1}a_i-a_ia_{s+1})W=0$,
	because $(a_{s+1}a_i-a_ia_{s+1})W\subseteq W_0$ for $i\le s$ and $a_jW_0=0$ for $j>s$.
	\item $(a_sa_ja_i-a_ja_sa_i)W=(a_sa_j-a_ja_s)a_iW=0$,
	because $a_iW\subseteq W_+$ for $i\le s$ and $(a_sa_j-a_ja_s)W_+=0$ for $j>s$.
	\item $(a_ja_ia_s-a_ja_sa_i)W=a_j(a_ia_s-a_sa_i)W=0$,
	because $(a_ia_s-a_sa_i)W\subseteq W_0$ for $i>s$ and $a_jW_0=0$ for $j>s$.
	\item $(a_{s+1}a_{s+1}a_i-a_{s+1}a_ia_{s+1})W=a_{s+1}(a_{s+1}a_i-a_ia_{s+1})W=0$,
	because $(a_{s+1}a_i-a_ia_{s+1})W\subseteq W_0$ for $i\le s$, and $a_{s+1}W_0=0$.
	\item $(a_sa_ia_s-a_ia_sa_s)W=(a_sa_i-a_ia_s)a_sW=0$, because $a_sW\subseteq W_+$
	and $(a_sa_i-a_ia_s)W_+=0$ for $i>s$.
\end{enumerate}

Now, let us fix $0\ne w\in W$ and $i_1,\dotsc,i_s\ge 0$. To show
that $xw=e_{i_1,\dotsc,i_s}$ for some $x\in K\cdot C_n$ write
$w=\sum_{j_1,\dotsc,j_s=0}^r\lambda_{j_1,\dotsc,j_s}e_{j_1,\dotsc,j_s}$,
where $\lambda_{j_1,\dotsc,j_s}\in K$. Define
\[m=\max\{j_s:\lambda_{j_1,\dotsc,j_s}\ne 0\text{ for some }j_1,\dotsc,j_{s-1}\}.\]
Then replacing the vector $w$ by $a_{s+1}^mw$ we may assume that
$\lambda_{j_1,\dotsc,j_s}=0$ for all $j_1,\dotsc,j_{s-1}\ge 0$ and $j_s>0$, that is
\[0\ne w=\sum_{j_1,\dotsc,j_{s-1}=0}^r\lambda_{j_1,\dotsc,j_{s-1},0}e_{j_1,\dotsc,j_{s-1},0}=u\otimes f_0,\]
where $0\ne u=\sum_{j_1,\dotsc,j_{s-1}=0}^r\lambda_{j_1,\dotsc,j_{s-1},0}e_{j_1,\dotsc,j_{s-1}}\in U$.
By assumptions on $U$, there exists $x\in K\cdot C_n^{s,s+1}$ such that
$xu=e_{i_1,\dotsc,i_{s-1}}$. Let $p=\sum_{j>s+1}\deg_{a_j}(x)$ and
$q=\sum_{j<s}\deg_{a_j}(x)$. Then
\begin{align*}
	\lambda_s^{-p-i_s}a_s^{i_s}a_{s+1}^qxa_s^pw & =(\lambda_s^{-1}a_s)^{i_s}a_{s+1}^q
	x(\lambda_s^{-1}a_s)^p(u\otimes f_0)\\
	& = (\lambda_s^{-1}a_s)^{i_s}a_{s+1}^qx(u\otimes f_p)\\
	& =(\lambda_s^{-1}a_s)^{i_s}a_{s+1}^qe_{i_1,\dotsc,i_{s-1},q}\\
	& =(\lambda_s^{-1}a_s)^{i_s}e_{i_1,\dotsc,i_{s-1},0}\\
	& =e_{i_1,\dotsc,i_s}.
\end{align*}
Hence the result follows.
\end{proof}
\end{proposition}

The next result is one of the essential tools used in this section.
An easy proof, based on the Density Theorem, can be found in \cite{cedo-okn}.

\begin{proposition}\label{central}
Let $A$ be a left primitive algebra over an algebraically closed field $F$.
If $\dim_FA<|F|$ then the algebra $A$ is central (that is, $Z(A)=F$).
\end{proposition}

If $w_1,w_2,\dotsc,w_k\in M$ for a monoid $M$ then we will write
$w_1^*w_2^*\dotsm w_k^*$ for the set of all elements of the form
$w_1^{i_1}w_2^{i_2}\dotsm w_k^{i_k}\in M$ with $i_1,i_2,\dotsc,i_k\ge 0$.

Now, consider a simple left $K[C_n]$-module $V$ with annihilator
$P$. Since $P$ is a prime ideal, it follows that $P$ contains a
minimal prime ideal of $K[C_n]$, which is of the form
$I_{\rho(d)}$ for some leaf $d$ in $D$ (see Section~\ref{background}
or \cite{chin2}). So, it is reasonable to investigate the
structure of left primitive ideals of $K[C_n]$ containing ideals
coming from diagrams of a particular shape. Our first result in
this direction reads as follows.

\begin{proposition}\label{prim1}
Assume that $P$ is a prime ideal of $K[C_n]$ containing the ideal
$I_\rho$, where $\rho$ is the congruence on $C_n$ determined by
the diagram

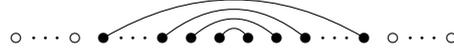
\[\diag{
	\circ & \dotsm & \circ & \bullet \ar@/^1.2pc/@{-}[rrrrrrrrr]
	& \dotsm & \bullet \ar@/^0.9pc/@{-}[rrrrr] & \bullet \ar@/^0.6pc/@{-}[rrr]
	& \bullet \ar@/^0.3pc/@{-}[r] & \bullet & \bullet & \bullet
	& \dotsm & \bullet & \circ & \dotsm & \circ
}\]
consisting of $t>0$ consecutive arcs $\arc{a_{s+1}a_s},\dotsc,\arc{a_{s+t}a_{s-t+1}}$
(as shown in the picture). Assume additionally that $a_{s+t}a_{s-t+1}\in P$.
Then $a_{s+t}\in P$ or $a_{s-t+1}\in P$.
\begin{proof}
Let $T$ be the image of $C_n$ in $K[C_n]/P$. 
Our first aim is to show that $a_{s-t+1}a_{s+t}=0$ in $T$. To prove this let us
introduce some notation. For any $1\le i\le j\le n$ let $W_{i,j}$
denote the subset of $C_n$ consisting of all elements of the form
$b_i\dotsm b_j$ written in the notation of the canonical
form~(\ref{canonical}). Moreover, let us adopt the convention that
$W_{i,j}=\{1\}$ in case $i>j$. In the following we shall use the
same notation for the elements of $T$.

First, using $a_{s+t}a_{s-t+1}=0$ in $T$, note that:
\begin{itemize}
	\item If $j\le s-t+1$ then $(a_{s-t+1}a_{s+t})a_j=a_{s+t}a_{s-t+1}a_j$ in $C_n$.
	Hence $(a_{s-t+1}a_{s+t})a_j=0$ in $T$ for all $j\le s-t+1$.
	\item If $j<k$ satisfy $j\le s-t+1$ and $k\le s+t$ then
	$(a_{s-t+1}a_{s+t})(a_ka_j)=a_{s-t+1}a_{s+t}a_ja_k=a_{s+t}a_{s-t+1}a_ja_k$ in $C_n$.
	Hence $(a_{s-t+1}a_{s+t})(a_ka_j)=0$ in $T$ for all $j<k$ such that $j\le s-t+1$ and $k\le s+t$.
\end{itemize}
This implies that $a_{s-t+1}a_{s+t}W_{1,s-t+1}=0$ in $T$ and
$a_{s-t+1}a_{s+t}W_{s-t+2,s+t}=a_{s-t+1}a_{s+t}U\cup\{0\}$ in $T$,
where $U=U_1\dotsm U_{2t-1}$ and
\begin{align*}
	U_1 & =(a_{s-t+2})^*,\\
	U_2 & =(a_{s-t+3}a_{s-t+2})^*(a_{s-t+3})^*,\\
	U_3 & =(a_{s-t+4}a_{s-t+2})^*(a_{s-t+4}a_{s-t+3})^*(a_{s-t+4})^*,\\
	& \vdotswithin{=}\\
	U_{2t-1} & =(a_{s+t}a_{s-t+2})^*(a_{s+t}a_{s-t+3})^*\dotsm(a_{s+t}a_{s+t-1})^*(a_{s+t})^*.
\end{align*}
As a consequence we get
\begin{equation}\label{prim1-1}
	a_{s-t+1}a_{s+t}T\subseteq a_{s-t+1}a_{s+t}UW_{s+t+1,n}\cup\{0\}.
\end{equation}
Next, remembering that $a_{s+t}a_{s-t+1}=0$ in $T$, we get:
\begin{itemize}
	\item If $j\ge s+t$ then $a_j(a_{s-t+1}a_{s+t})=a_ja_{s+t}a_{s-t+1}$ in $C_n$.
	Hence $a_j(a_{s-t+1}a_{s+t})=0$ in $T$ for all $j\ge s+t$.
	\item If $j<k$ satisfy $j\le s-t+1$ and $k\ge s+t$ then
	$(a_ka_j)(a_{s-t+1}a_{s+t})=(a_{s-t+1}a_{s+t})(a_ka_j)$ in $C_n$.
	Hence also $(a_ka_j)(a_{s-t+1}a_{s+t})=(a_{s-t+1}a_{s+t})(a_ka_j)$
	in $T$ for all $j<k$ such that $j\le s-t+1$ and $k\ge s+t$.
	\item If $j<k$ satisfy $j>s-t+1$ and $k\ge s+t$ then
	$(a_ka_j)(a_{s-t+1}a_{s+t})=(a_ka_{s-t+1})(a_ja_{s+t})=(a_ja_{s+t})(a_ka_{s-t+1})=a_ja_ka_{s+t}a_{s-t+1}$
	in $C_n$. Hence $(a_ka_j)(a_{s-t+1}a_{s+t})=0$ in $T$ for all $j<k$ such that $j>s-t+1$ and $k\ge s+t$.
\end{itemize}
These equalities assure that
$W_{s+t+1,n}a_{s-t+1}a_{s+t}\subseteq a_{s-t+1}a_{s+t}T$.
Thus, together with (\ref{prim1-1}), we obtain
\begin{equation}
	a_{s-t+1}a_{s+t}Ta_{s-t+1}a_{s+t}
	\subseteq a_{s-t+1}a_{s+t}UW_{s+t+1,n}a_{s-t+1}a_{s+t}
	\subseteq a_{s-t+1}a_{s+t} Ua_{s-t+1}a_{s+t}T.\label{prim1-2}
\end{equation}
Now, choose $1\ne u\in U$. Then let $m$ be the minimum of those
numbers $j\in\{s-t+2,\dotsc,s+t\}$ such that the generator $a_j$
appears in $u$. Since $a_{s+t}(a_ka_m)=a_k(a_{s+t}a_m)$ in $C_n$ for all
$m<k\le s+t$, we get $a_{s+t}u\in Ta_{s+t}a_m$. Therefore,
\begin{equation}
	a_{s-t+1}a_{s+t}ua_{s-t+1}a_{s+t}\in a_{s-t+1}T(a_{s+t}a_ma_{s-t+1})a_{s+t}
	=a_{s-t+1}T(a_{s+t}a_{s-t+1}a_m)a_{s+t}=0.\label{prim1-3}
\end{equation}
Thus (\ref{prim1-3}) together with $(a_{s-t+1}a_{s+t})^2=0$ in $T$ yield
$a_{s-t+1}a_{s+t}Ua_{s-t+1}a_{s+t}=0$ in $T$. Hence, as a consequence of (\ref{prim1-2}),
we get $a_{s-t+1}a_{s+t}Ta_{s-t+1}a_{s+t}=0$ as well. Since $T$ is a
prime semigroup, we conclude that $a_{s-t+1}a_{s+t}=0$ in $T$, as desired.

Next, we claim that $a_{s-t+1}Ta_{s+t}=0$. First, by (\ref{rel2-1})
it follows that $a_{s-t+1}$ commutes in $T$ with $a_1,\dotsc,a_s$.
Hence we get $a_{s-t+1}W_{1,s}=W_{1,s}a_{s-t+1}$ in $T$.
Similarly, by (\ref{rel2-2}) it follows that $a_{s+t}$ commutes in $T$
with $a_{s+1},\dotsc,a_n$. Moreover, $a_{s+t}$ commutes in $C_n$
with $a_ja_i$ for all $i<j$ such that $i\le s+t\le j$. Therefore, we get
$W_{s+t,n}a_{s+t}=a_{s+t}W_{s+t,n}$ in $T$, which leads to
\[a_{s-t+1}Ta_{s+t}=W_{1,s}a_{s-t+1}W_{s+1,n}a_{s+t}
=W_{1,s}a_{s-t+1}W_{s+1,s+t-1}a_{s+t}W_{s+t,n}.\]
So it is enough to show that $a_{s-t+1}W_{s+1,s+t-1}a_{s+t}=0$ in $T$.
Further, let us observe that:
\begin{itemize}
	\item If $j>s$ then $a_ja_{s+t}=a_{s+t}a_j$ in $T$ (by (\ref{rel2-2}),
	because $j>s$ and $s+t>s$),
	\item If $k>s$ then $(a_ka_s)a_{s+t}=a_{s+t}a_sa_k$ in $T$
	(by (\ref{rel2-4}) with $r=1$, because $k>s$ and $s+t>s$),
	\item If $j<k$ satisfy $s-t<j\le s$ and $j+k>n$ then
	$(a_ka_j)a_{s+t}=a_{s+t}a_ja_k$ in $T$ (by (\ref{rel2-4}) with $r=s-j+1$.
	Indeed, the assumption $s-t<j\le s$ assures that $1\le r\le t$.
	Moreover, $s+r=2s-j+1=n-j+1$. Hence, to use (\ref{rel2-4}), it only
	remains to check that $k\ge n-j+1$ and $s+t\ge n-j+1$.
	Now, the first inequality is a consequence of $j+k>n$,
	whereas the second one is obtained as follows.
	Since $s-t<j$, we get $t\ge s-j+1$. Therefore, $s+t\ge 2s-j+1=n-j+1$. 
	We note that straightforward calculations on indices of this type will be also
	used in other proofs in this section; however, complete explanations will be skipped.)
\end{itemize}
These equalities lead to the conclusion that
$W_{s+1,s+t-1}a_{s+t}\subseteq Va_{s+t}T$,
where $V=V_1\dotsm V_{t-1}$ and
\begin{align*}
	V_1 & =(a_{s+1}a_1)^*(a_{s+1}a_2)^*\dotsm(a_{s+1}a_{s-1})^*,\\
	V_2 & =(a_{s+2}a_1)^*(a_{s+2}a_2)^*\dotsm(a_{s+2}a_{s-2})^*,\\
	V_3 & =(a_{s+3}a_1)^*(a_{s+3}a_2)^*\dotsm(a_{s+3}a_{s-3})^*,\\
	& \vdotswithin{=}\\
	V_{t-1} & =(a_{s+t-1}a_1)^*(a_{s+t-1}a_2)^*\dotsm(a_{s+t-1}a_{s-t+1})^*.
\end{align*}
Next, we have $a_{s-t+1}(a_ka_j)=a_ja_ka_{s-t+1}$ in $T$ for all
$j<k$ such that $j+k\le n+1$ and $s<k\le s+t$. Indeed, here
assumptions on $j$ and $k$ can be rewritten as $j\le n-k+1$ and
$s-t+1\le n-k+1$, hence our equality follows by (\ref{rel2-3}) with $r=k-s$.
Since each element of $V$ is a product of elements $a_ka_j$ with $j<k$
such that $j+k\le n$ and $s<k<s+t$, we conclude that
$a_{s-t+1}V\subseteq Ta_{s-t+1}$. Therefore
\[a_{s-t+1}W_{s+1,s+t-1}a_{s+t}\subseteq
a_{s-t+1}Va_{s+t}T\subseteq Ta_{s-t+1}a_{s+t}T=0.\]
This proves the claim. Since $T$ is a  prime semigroup,
it follows that $a_{s-t+1}=0$ in $T$ or $a_{s+t}=0$ in $T$
or, in other words, $a_{s-t+1}\in P$ or $a_{s+t}\in P$.
\end{proof}
\end{proposition}

\begin{proposition}\label{prim2}
Assume that $P$ is a left primitive ideal of $K[C_n]$ containing
the ideal $I_\rho$, where $\rho$ is the congruence on $C_n$
determined by one of the two diagrams

\[\diag{
	\circ & \dotsm & \circ & \bullet & \bullet \ar@/^0.9pc/@{-}[rrrrrrr]
	& \dotsm & \bullet \ar@/^0.6pc/@{-}[rrr] & \bullet \ar@/^0.3pc/@{-}[r]
	& \bullet & \bullet & \dotsm & \bullet & \circ & \dotsm & \circ}
	\qquad\text{or}\qquad\diag{\circ & \dotsm & \circ & \bullet
	\ar@/^0.9pc/@{-}[rrrrrrr] & \dotsm & \bullet \ar@/^0.6pc/@{-}[rrr]
	& \bullet \ar@/^0.3pc/@{-}[r] & \bullet & \bullet & \dotsm
	& \bullet & \bullet & \circ & \dotsm & \circ
}\]
consisting of $t>0$ consecutive arcs
$\arc{a_{s+1}a_s},\dotsc,\arc{a_{s+t}a_{s-t+1}}$ and a single dot $a_{s-t}$
or $a_{s+t+1}$ (as shown in the picture). Assume additionally that
$a_{s+t}a_{s-t+1}\notin P$. Then $a_{s-t}-\lambda a_{s-t+1}\in P$ or
$a_{s+t+1}-\lambda a_{s+t}\in P$ for some $\lambda\in K$, respectively.
\begin{proof}
Since the two cases are symmetric, it is enough to
consider the former. Let us notice first that the element
$a_{s+t}a_{s-t+1}$ is central in $K[C_n]/P$. Indeed, if $i\le s-t+1$,
then $a_ia_{s+t}a_{s-t+1}=a_{s-t+1}a_{s+t}a_i=a_{s+t}a_{s-t+1}a_i$ in $K[C_n]/P$
(first equality is a consequence of (\ref{rel3-4}) with $r=t$, because $i\le s-t+1$;
second equality is valid in $C_n$). Next, if $s-t+1<i\le s+t$, then
$a_ia_{s+t}a_{s-t+1}=a_{s+t}a_{s-t+1}a_i$ in $C_n$, hence also in $K[C_n]/P$.
Whereas, if $i>s+t$, then $a_ia_{s+t}a_{s-t+1}=a_ia_{s-t+1}a_{s+t}=a_{s+t}a_{s-t+1}a_i$
in $K[C_n]/P$ (first equality holds in $C_n$; second equality follows from (\ref{rel3-5})
with $r=t$, because $i\ge s+t$). Therefore, Proposition~\ref{central}
implies that $a_{s+t}a_{s-t+1}=\mu$ in $K[C_n]/P$ for some $\mu\in K$.
Moreover, $\mu\ne 0$ because $a_{s+t}a_{s-t+1}\notin P$. We
claim that the element $a_{s+t}a_{s-t}$ is central in $K[C_n]/P$
as well. Indeed, if $j\le s-t$, then
$a_ja_{s+t}a_{s-t}=a_{s-t}a_{s+t}a_j=a_{s+t}a_{s-t}a_j$ in $K[C_n]/P$
(first equality holds by (\ref{rel3-4}) with $r=t$, because $j\le s-t+1$ and
$s-t\le s-t+1$; second equality is valid in $C_n$). Next, if $s-t<j<s+t$, then
$a_ja_{s+t}a_{s-t}=a_{s+t}a_{s-t}a_j$ in $C_n$, hence also in $K[C_n]/P$.
Whereas, if $j\ge s+t$, then
$a_ja_{s+t}a_{s-t}=a_{s+t}a_ja_{s-t}=a_{s+t}a_{s-t}a_j$ in $K[C_n]/P$
(first equality holds $C_n$; second equality is a consequence of
(\ref{rel3-3}), because $j\ge s+t$). Therefore, by Proposition~\ref{central},
we get $a_{s+t}a_{s-t}=\nu$ in $K[C_n]/P$ for some $\nu\in K$, which leads to
\[\mu a_{s-t}=a_{s+t}a_{s-t+1}a_{s-t}=a_{s+t}a_{s-t}a_{s-t+1}=\nu a_{s-t+1}\]
in $K[C_n]/P$. Hence we get $a_{s-t}-\lambda a_{s-t+1}\in P$, where
$\lambda=\mu^{-1}\nu\in K$.
\end{proof}
\end{proposition}

\begin{proposition}\label{prim3}
Assume that $P$ is a left primitive ideal of $K[C_n]$ containing
the ideal $I_\rho$, where $\rho$ is the congruence on $C_n$
determined by one of the two diagrams

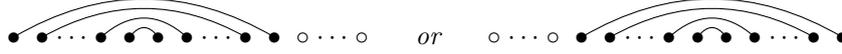
\[\diag{
	\bullet \ar@/^1.2pc/@{-}[rrrrrrrrr] & \bullet \ar@/^0.9pc/@{-}[rrrrrrr]
	& \dotsm & \bullet \ar@/^0.6pc/@{-}[rrr] & \bullet \ar@/^0.3pc/@{-}[r]
	& \bullet & \bullet & \dotsm & \bullet & \bullet & \circ & \dotsm & \circ}
	\qquad\text{or}\qquad\diag{ \circ & \dotsm & \circ & \bullet
	\ar@/^1.2pc/@{-}[rrrrrrrrr] & \bullet \ar@/^0.9pc/@{-}[rrrrrrr]
	& \dotsm & \bullet \ar@/^0.6pc/@{-}[rrr] & \bullet \ar@/^0.3pc/@{-}[r]
	& \bullet & \bullet & \dotsm & \bullet & \bullet
}\]
consisting of $t>0$ consecutive arcs
$\arc{a_{t+1}a_t},\dotsc,\arc{a_{2t}a_1}$ or
$\arc{a_{n-t+1}a_{n-t}},\dotsc,\arc{a_na_{n-2t+1}}$
(as shown in the picture). Assume additionally that $a_{2t}a_1\notin P$
or $a_na_{n-2t+1}\notin P$, respectively. Then
$a_{2t+1}-\lambda a_{2t}\in P$ or $a_{n-2t}-\lambda a_{n-2t+1}\in P$
for some $\lambda\in K$, respectively.
\begin{proof}
Since the two cases are symmetric, it suffices to consider
the former. First, notice that the element $a_{2t}a_1$ is central
in $K[C_n]/P$. Indeed, if $i\le 2t$, then
$a_ia_{2t}a_1=a_{2t}a_1a_i$ in $C_n$, hence also in $K[C_n]/P$.
Whereas, if $i>2t$, then $a_ia_{2t}a_1=a_ia_1a_{2t}=a_{2t}a_1a_i$
(first equality holds in $C_n$; second equality follows from (\ref{rel2-4})
with $r=s=t$, because $i\ge 2t$). Hence, by Proposition~\ref{central},
we get $a_{2t}a_1=\mu$ in $K[C_n]/P$ for some $\mu\in K$. Moreover,
$a_{2t}a_1\notin P$ implies that $\mu\ne 0$. We claim that the
element $a_{2t+1}a_1$ is central in $K[C_n]/P$ as well. Indeed, if
$j\le 2t$, then $a_{2t+1}a_1a_j=a_ja_{2t+1}a_1$ in $C_n$,
hence also in $K[C_n]/P$. Whereas, if $j>2t$, then
$a_ja_{2t+1}a_1=a_ja_1a_{2t+1}=a_{2t+1}a_1a_j$ in $K[C_n]/P$
(first equality holds in $C_n$; second equality is a consequence of
(\ref{rel2-4}) with $r=s=t$, because $j\ge 2t$).
Thus, Proposition~\ref{central} yields $a_{2t+1}a_1=\nu$ in
$K[C_n]/P$ for some $\nu\in K$, and we get
\[\mu a_{2t+1}=a_{2t+1}a_{2t}a_1=a_{2t+1}a_1a_{2t}=\nu a_{2t}\]
in $K[C_n]/P$. Hence we get $a_{2t+1}-\lambda a_{2t}\in P$,
where $\lambda=\mu^{-1}\nu\in K$.
\end{proof}
\end{proposition}

Before we proceed to a formulation of the main result of this
paper, let us recall some notions and introduce some notation.
We say that the element $x\in K[C_n]$ acts regularly on a left
$K[C_n]$-module $V$ if $xv\ne 0$ for each $0\ne v\in V$ (or, in
other words, if the annihilator of $x$ in $V$ is equal to zero).

In the following five lemmas we assume that $n$ is even, say
$n=2s$ for some $s\ge 1$. We assume as well that the
$K[C_n]$-module $V$ is simple, and its annihilator $P$ contains
the ideal $I_\rho$, where $\rho$ is the  congruence on $C_n$
determined by the diagram\medskip

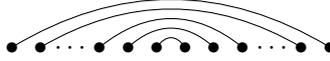
\[\diag{
	\bullet \ar@/^1.5pc/@{-}[rrrrrrrrrrr] & \bullet \ar@/^1.2pc/@{-}[rrrrrrrrr]
	& \dotsm & \bullet \ar@/^0.9pc/@{-}[rrrrr] & \bullet \ar@/^0.6pc/@{-}[rrr]
	& \bullet \ar@/^0.3pc/@{-}[r] & \bullet & \bullet & \bullet & \dotsm & \bullet & \bullet
}\]
consisting of $s$ consecutive arcs $\arc{a_{s+1}a_s},\dotsc,
\arc{a_na_1}$ (as shown in the picture). Moreover, we consider the set
\[X=\{a_ja_i:i<s+1<j\text{ and }i+j>n+1\}\subseteq C_n,\]
and its subsets
\[X_0=\{x\in X:x\text{ does not act regularly on }V\},\qquad X_1=X\setminus X_0.\]
Of course, the sets $X_0$ and $X_1$ depend on the module $V$.

\begin{lemma}\label{lem1}
We have $xy-yx\in P$ for all $x,y\in X$.
\begin{proof}
Let $x=a_ja_i$ and $y=a_la_k$. If $j=l$ then $xy=yx$ in $C_n$.
Hence, we may assume that $j<l$. If $i\ge k$ then again $xy=yx$ in
$C_n$. So, assume that $i<k$. Summarizing, we are in the
situation where $i<k<s+1<j<l$. Then, by (\ref{rel2-4}) with $r=s-i+1$,
we get $a_ja_ia_l=a_la_ia_j$ in $K[C_n]/P$ (because $l>j>n-i+1$
so, in particular, $j,l\ge n-i+1$), which leads to
\[xy=a_ja_ia_la_k=a_la_ia_ja_k=a_la_ja_ia_k=a_la_ka_ja_i=yx\]
in $K[C_n]/P$. Hence the result follows.
\end{proof}
\end{lemma}

\begin{lemma}\label{lem2}
Assume that $a_{n-j+1}a_j\notin P$ for each $j=1,\dotsc,s$.
If $x\in X_1$ then $x-\mu\in P$ for some $0\ne\mu\in K$.
\begin{proof}
First, notice that each element $a_{n-j+1}a_j$ for $j=1,\dotsc,s$
is central in $K[C_n]/P$. Indeed, if $i\le j$, then
$a_ia_{n-j+1}a_j=a_ja_{n-j+1}a_i=a_{n-j+1}a_ja_i$ in $K[C_n]/P$
(first equality is a consequence of (\ref{rel2-3}) with $r=s-j+1$, because $i\le j$;
second equality is valid in $C_n$). Next, if $j<i\le n-j+1$, then $a_i$
commutes with $a_{n-j+1}a_j$ in $C_n$, hence also in $K[C_n]/P$.
Whereas, if $i>n-j+1$, then
$a_ia_{n-j+1}a_j=a_ia_ja_{n-j+1}=a_{n-j+1}a_ja_i$ in $K[C_n]/P$
(first equality holds in $C_n$; second equality follows from (\ref{rel2-4})
with $r=s-j+1$, because $i\ge n-j+1$).
Hence Proposition~\ref{central} assures that $a_{n-j+1}a_j=\lambda_j$
in $K[C_n]/P$. Since $a_{n-j+1}a_j\notin P$, we get $\lambda_j\ne 0$
for each $j=1,\dotsc,s$. Now, let $x=a_ka_j$. Since we have
\[(a_ka_j)(a_{n-j+1}a_{n-k+1})=a_{n-j+1}a_ka_ja_{n-k+1}=(a_{n-j+1}a_j)(a_ka_{n-k+1})\]
in $C_n$ and $a_ka_j(a_{n-j+1}a_{n-k+1}a_ka_j-\lambda_j\lambda_{n-k+1})V=0$,
we conclude that $x$ and $y=\nu a_{n-j+1}a_{n-k+1}$, where
$\nu=(\lambda_j\lambda_{n-k+1})^{-1}\in K$,
are mutual inverses in $K[C_n]/P$.

We claim that each generator $a_i$, for $i=1,\dotsc,n$,
commutes with $x$ or $y$ in $K[C_n]/P$. Indeed, if $i\le j$
then, by (\ref{rel2-3}) with $r=s-j+1$, we get
$a_ia_{n-j+1}a_{n-k+1}=a_{n-k+1}a_{n-j+1}a_i$ in $K[C_n]/P$
(because we have $i\le j$, and from $j+k>n+1$ we get also $n-k+1\le j$), so
\[a_iy=\nu a_ia_{n-j+1}a_{n-k+1}=\nu a_{n-k+1}a_{n-j+1}a_i=\nu a_{n-j+1}a_{n-k+1}a_i=ya_i\]
in $K[C_n]/P$. If $j<i\le n-j+1$ then from $j+k>n+1$ we get $n-k+1<j<i$,
so $a_i$ commutes with $y=a_{n-j+1}a_{n-k+1}$ in $C_n$,
hence in $K[C_n]/P$, too. If $n-j+1<i\le k$ then $j<n-j+1<i$ implies that 
$a_i$ commutes with $x=a_ka_j$ in $C_n$, hence in $K[C_n]/P$ as well.
Whereas, if $i>k$ then, by (\ref{rel2-4}) with $r=s-j+1$, we have
$a_ia_ja_k=a_ka_ja_i$ in $K[C_n]/P$ (because from $j+k>n+1$
we get $i>k\ge n-j+1$), which yields
\[a_ix=a_ia_ka_j=a_ia_ja_k=a_ka_ja_i=xa_i\]
in $K[C_n]/P$, and the claim follows.
In particular, $x$ in central in $K[C_n]/P$. Therefore,
Proposition~\ref{central} guarantees that $x-\mu\in P$ for some $\mu\in K$.
Since $x$ acts regularly on $V$, we must have $\mu\ne 0$.
\end{proof}
\end{lemma}

Let $x=a_la_k\in X$ and $y=a_ja_i\in X$. We say that $x$ dominates
$y$ (or that $y$ is dominated by $x$) if $j\le l$ and $i\le k$.

\begin{lemma}\label{lem3}
Assume that $a_{n-j+1}a_j\notin P$ for each $j=1,\dotsc,s$.
If $y\in X$ is dominated by some $x\in X_1$, then $y\in X_1$.
\begin{proof}
By Lemma~\ref{lem2} we know that $x-\mu\in P$ for some $0\ne\mu\in
K$. Moreover, as in Lemma~\ref{lem2}, we have $a_{n-j+1}a_j=\lambda_j$
in $K[C_n]/P$, where $\lambda_j\ne 0$, for each $j=1,\dotsc,s$.
Now, write $x=a_la_k$ and $y=a_ja_i$. Of
course it suffices to consider just two cases. Namely,
$(i,j)=(k-1,l)$ and $(i,j)=(k,l-1)$. Assume first that
$(i,j)=(k-1,l)$. Then, by (\ref{rel2-4}) with $r=s-k+1$, we have
$a_{n-k+2}a_ka_l=a_la_ka_{n-k+2}$ in $K[C_n]/P$ (because from $k+l>n+1$
we get $l\ge n-k+1$ and, of course, $n-k+2\ge n-k+1$). Now, if $v\in V$ and $yv=0$ then
\[0=a_{n-k+2}a_kyv=a_{n-k+2}a_ka_la_{k-1}v=a_la_ka_{n-k+2}a_{k-1}v=\lambda_{k-1}\mu v,\]
because $a_la_k=\mu$ and $a_{n-k+2}a_{k-1}=\lambda_{k-1}$ in $K[C_n]/P$.
Hence $v=0$. Finally, let $(i,j)=(k,l-1)$. If $yv=0$ for some $v\in V$ then
\[0=a_la_{n-l+2}yv=a_la_{n-l+2}a_{l-1}a_kv=a_la_{l-1}a_{n-l+2}a_kv
=\lambda_{n-l+2}a_la_kv=\lambda_{n-l+2}\mu v,\]
because $a_{l-1}a_{n-l+2}=\lambda_{n-l+2}$ and $a_la_k=\mu$ in $K[C_n]/P$.
Hence again $v=0$, and the result follows.
\end{proof}
\end{lemma}

The statement of the last lemma can be easily visualized if we
arrange the elements of $X$ in a triangular matrix, as follows
\begin{equation}\begin{array}{lllllll}\label{matrix}
	a_na_2 & a_na_3        & a_na_     4  & \dotsm & a_na_{s-2}       & a_na_{s-1}      & a_na_s\\
	& a_{n-1}a_3  & a_{n-1}a_4 & \dotsm & a_{n-1}a_{s-2}  & a_{n-1}a_{s-1} & a_{n-1}a_s\\
	&                     & a_{n-2}a_4 & \dotsm & a_{n-2}a_{s-2}  & a_{n-2}a_{s-1} & a_{n-2}a_s\\
	& & & \ddots  & \multicolumn{1}{c}{\vdots} & \multicolumn{1}{c}{\vdots} & \multicolumn{1}{c}{\vdots}\\
	&                     &                    &             & a_{s+4}a_{s-2} & a_{s+4}a_{s-1} & a_{s+4}a_s\\
	&                     &                    &            &                           & a_{s+3}a_{s-1} & a_{s+3}a_s\\
	&                      &                   &            &                           &                           & a_{s+2}a_s
\end{array}\end{equation}
Then, for each $x\in X$, elements in $X$ dominated by $x$
constitute a right triangle with $x$ as the vertex of the right
angle and with its hypotenuse consisting of elements lying on the
diagonal of the above matrix. It is also worth to reformulate
Lemma~\ref{lem3} in the following way. If $x=a_ja_i\in X_0$ then
all elements dominating $x$ also lie in $X_0$ (these are precisely
the elements of the matrix lying inside the rectangle defined by
the vertices $a_ja_i$, $a_na_i$, $a_na_s$, $a_ja_s$).

\begin{lemma}\label{lem4}
    If $x\in X_0$ then for each $v\in V$ there exists $m>0$ such that $x^mv=0$.
    \begin{proof}
Let $x=a_ka_j$. Since $x$ does not act regularly on $V$, there
exists $0\ne w\in V$ such that $xw=0$. Because $V$ is a simple
module, we have $V=K[C_n]w$. We claim that for each $z\in C_n$
there exists $l>0$ such that $x$ commutes with $x^lz$ in
$K[C_n]/P$ (in fact, it suffices to take $l=\deg(z)$). Of course
this claim implies our lemma. We shall prove the claim by
induction on $\deg(z)$. So assume first that $\deg(z)=1$. Then
$z=a_i$ for some $i\in\{1,\dotsc,n\}$. If $i<j$ then
$xz=a_ka_ja_i=(a_ka_i)a_j$ in $C_n$. Hence $xz$ commutes with $x$ in $C_n$,
so also in $K[C_n]/P$. Next, if $j\le i\le k$ then $x$ commutes with $z$ in $C_n$,
hence also with $xz$ in $C_n$, and of course in $K[C_n]/P$ as well. Finally,
if $i>k$ then, by (\ref{rel2-4}) with $r=s-j+1$, we get $xz=a_ka_ja_i=(a_ia_j)a_k$
in $K[C_n]/P$ (because $j+k>n+1$ implies $i>k\ge n-j+1$).
Thus $x$ commutes with $xz$ in $K[C_n]/P$ as well. Now
assume that $\deg(z)>1$ and write $z=z'a_i$ for some $z'\in C_n$
and some $i\in\{1,\dotsc,n\}$. By induction, there exists $l>0$
such that $x$ commutes with $x^lz'$ in $K[C_n]/P$. Now
$x^{l+1}z=x^{l+1}z'a_i=(x^lz')(xa_i)$, so the claim follows,
because $x$ commutes in $K[C_n]/P$ with $x^lz'$ and with $xa_i$,
hence also with $x^{l+1}z$.
\end{proof}
\end{lemma}

\begin{lemma}\label{lem5}
Assume that $a_{n-j+1}a_j\notin P$ for each $j=1,\dotsc,s$.
If $a_{s+1}$ does not act regularly on $V$ then there exists
$0\ne v\in V$ such that $xv=0$ for each $x\in X_0$ and $a_jv=0$
for each $j>s$. In this situation $V$ is spanned as a $K$-linear space by
the set
\[(a_{n-1}a_1)^*(a_{n-2}a_2)^*(a_{n-3}a_3)^*\dotsm
(a_{s+3}a_{s-3})^*(a_{s+2}a_{s-2})^*(a_{s+1}a_{s-1})^*a_s^*v.\]
\begin{proof}
First, note that for each $j=1,\dotsc,s$ we have
$a_{n-j+1}a_j=\lambda_j$ in $K[C_n]/P$ for some $0\ne\lambda_j\in K$
(this was already proved in Lemma~\ref{lem2}). Assume for the moment
that we already have a vector $0\ne w\in V$ such that $xw=0$ for each
$x\in X_0$. Then there exists $k>0$ such that
$a_{s+1}^kw=0\ne a_{s+1}^{k-1}w$
(the proof of this fact is completely analogous to
the proof of Lemma~\ref{lem4}, so it will be omitted here). Let
$v=a_{s+1}^{k-1}w\ne 0$. Since $a_{s+1}$ commutes with each
$x\in X$ in $C_n$, we get $xv=0$ for each $x\in X_0$. Of course
$a_{s+1}v=0$, and if $j>s+1$ then we have
\[0=a_ja_sa_{s+1}v=a_ja_{s+1}a_sv=\lambda_sa_jv,\]
which gives $a_jv=0$, because $\lambda_s\ne 0$. Thus, to finish
the proof of the first part of our lemma, it is enough to show
that there exists $0\ne w\in V$ such that $xw=0$ for each $x\in X_0$.
If $X_0$ is empty then there is nothing to show. Therefore,
assume that $X_0=\{x_1,\dotsc,x_d\}$ with $d=|X_0|>0$. Take $l<d$
and suppose that there exists $0\ne w_l\in V$ such that
$x_1w_l=\dotsb=x_lw_l=0$. By Lemma~\ref{lem4} we know that
$x_{l+1}^mw_l=0\ne x_{l+1}^{m-1}w_l$ for some $m>0$. Then define
$w_{l+1}=x_{l+1}^{m-1}w_l\ne 0$. Because $x_1,\dotsc,x_l$ commute
with $x_{l+1}$ in $K[C_n]/P$ (see Lemma~\ref{lem1}), we get
$x_1w_{l+1}=\dotsb=x_lw_{l+1}=x_{l+1}w_{l+1}=0$. Now, it is clear
that after $d$ steps we obtain a nonzero vector $w=w_d\in V$ such
that $xw=0$ for each $x\in X_0$.

Let us proceed to the proof of the last statement of the lemma.
So, fix $0\ne v\in V$ satisfying $xv=0$ for each $x\in X_0$ and
$a_jv=0$ for each $j>s$. Of course $V$ is spanned as a $K$-linear
space by the set $C_nv$. Hence it suffices to show that for each
$x=b_1\dotsm b_n\in C_n$ written in its canonical form (\ref{canonical})
we have $xv\in K\cdot(a_{n-1}a_1)^*\dotsm(a_{s+1}a_{s-1})^*a_s^*v$.

First, by (\ref{rel2-1}), we have $a_ia_j=a_ja_i$ in
$K[C_n]/P$ for all $i,j\le s$. Hence, the element
$b_1\dotsm b_s$ can be written in $K[C_n]/P$ as an element of the
set $a_1^*\dotsm a_s^*$. Next, for each $j<s$ we have
\[\lambda_{j+1}a_j=a_{n-j}a_{j+1}a_j=(a_{n-j}a_j)a_{j+1}\]
in $K[C_n]/P$. Thus, we conclude that
$a_1^*\dotsm a_s^*\subseteq K\cdot(a_{n-1}a_1)^*\dotsm(a_{s+1}a_{s-1})^*a_s^*$
in $K[C_n]/P$, which allows us to assume that
\begin{equation}
	b_1\dotsm b_s\in (a_{n-1}a_1)^*\dotsm(a_{s+1}a_{s-1})^*a_s^*.\label{b1}
\end{equation}
Further, by (\ref{rel2-2}), we have $a_ka_l=a_la_k$
in $K[C_n]/P$ for all $k,l>s$. Therefore, for each $j>s$, the element $b_j$ can
be written in $K[C_n]/P$ as an element of the set
$(a_ja_1)^*\dotsm(a_ja_s)^*a_{s+1}^*\dotsm a_j^*$.
Because the elements
$a_{s+1},\dotsc,a_j$ commute in $C_n$ with each $a_ka_i$, where
$i<k$ satisfy $i\le s$ and $k>j$, we deduce that $a_{s+1},\dotsc,a_j$
commute in $K[C_n]/P$ with all elements $b_{j+1},\dotsc,b_n$.
Moreover, $a_{s+1}v=\dotsb=a_jv=0$. These two facts allow us to
assume that
\begin{equation}
	b_j\in(a_ja_1)^*\dotsm(a_ja_s)^*\text{ for each }j>s.\label{b2}
\end{equation}
Next, we claim that for all $i<s<j$ such that $i+j<n$ the equality
\begin{equation}\label{lem5-claim}
	(\lambda_{i+1}\dotsm\lambda_{n-j})a_ja_i=(a_{n-i}a_i)
	(a_{n-i-1}a_{i+1})(a_{n-i-2}a_{i+2})\dotsm(a_ja_{n-j})
\end{equation}
holds in $K[C_n]/P$.
We shall prove the claim by induction on $d=n-i-j$.
If $d=1$ then $i+j=n-1$ and we have
\[(a_{n-i}a_i)(a_ja_{n-j})=(a_{j+1}a_i)(a_ja_{i+1})
=a_ja_{j+1}a_ia_{i+1}=a_j(a_{j+1}a_{i+1})a_i=\lambda_{i+1}a_ja_i,\]
because $a_{j+1}a_{i+1}=a_{n-i}a_{i+1}=\lambda_{i+1}$ in $K[C_n]/P$.
So assume that $d>1$, and the claim is true for all $i<s<j$ such that
$n-i-j=d$. Our aim is to show that (\ref{lem5-claim}) holds for all $i<s<j$
such that $n-i-j=d+1$. Observe that in this case we must have $i+1<s$, because otherwise
$i\ge s-1$ and $j\ge s+1$ give $d+1=n-i-j\le n-(s-1)-(s+1)=0$, a contradiction.
So $i+1<s<j$ and $n-(i+1)-j=d$. Therefore, by induction, we get
\[(\lambda_{i+2}\dotsm\lambda_{n-j})a_ja_{i+1}=(a_{n-i-1}a_{i+1})\dotsm(a_ja_{n-j})\]
in $K[C_n]/P$. This equality, together with $a_{n-i}a_{i+1}=\lambda_{i+1}$ in $K[C_n]/P$,
lead to
\begin{align*}
	(\lambda_{i+1}\lambda_{i+2}\dotsm\lambda_{n-j})a_ja_i
	& =(\lambda_{i+2}\dotsm\lambda_{n-j})(a_{n-i}a_{i+1})(a_ja_i)\\
	& =(\lambda_{i+2}\dotsm\lambda_{n-j})a_ja_{n-i}a_{i+1}a_i\\
	& =(\lambda_{i+2}\dotsm\lambda_{n-j})a_ja_{n-i}a_ia_{j+1}\\
	& =(\lambda_{i+2}\dotsm\lambda_{n-j})(a_{n-i}a_i)(a_ja_{i+1})\\
	& =(a_{n-i}a_i)(a_{n-i-1}a_{i+1})\dotsm(a_ja_{n-j})
\end{align*}
in $K[C_n]/P$, hence the claim follows.
Now, looking at the form (\ref{b2}) of $b_j$ and using (\ref{lem5-claim})
to rewrite the factors $a_ja_i$ with $i<n-j$ (appearing in the form (\ref{b2}) of $b_j$),
and also noticing that $a_ja_{n-j+1}=\lambda_{n-j+1}$ in $K[C_n]/P$,
we conclude that we may restrict to the situation when
\begin{equation}
	b_j\in(a_{n-1}a_1)^*\dotsm(a_ja_{n-j})^*\cdot
	(a_ja_{n-j+2})^*\dotsm(a_ja_s)^*\text{ for each }j>s.\label{b3}
\end{equation}
Note that in the form (\ref{b3}) of $b_j$ two types of factors appear.
First $n-j$ factors are of the form $a_{n-i}a_i$ for $i=1,\dotsc,n-j$,
whereas next $s-(n-j+1)=j-s-1$ factors (separated from the first $n-j$ factors by a dot)
are of the form $a_ja_i$, where $i=n-j+2,\dotsc,s$.
Further, each factor $a_ja_i$ for $n-j+2\le i\le s$,
appearing in the form (\ref{b3}) of $b_j$, lies
in $X$ and commutes with all factors $a_{n-1}a_1,\dotsc,a_{j+1}a_{n-j-1}$
that appear in the form (\ref{b3}) of the elements $b_{j+1},\dotsc,b_n$.
Hence, we can write $b_{s+1}\dotsm b_n$ in $K[C_n]/P$
as an element of the set
$K\cdot(a_{n-1}a_1)^*\dotsm(a_{s+1}a_{s-1})^*\langle X\rangle$,
where $\langle X\rangle\subseteq C_n$ denotes the monoid generated
by the set $X=X_0\cup X_1$. Since $X_0v=0$ and $X_1v\subseteq Kv$ (see Lemma~\ref{lem2}),
we get $\langle X\rangle v\subseteq Kv$, which leads to the conclusion that
\begin{equation}
	b_{s+1}\dotsm b_nv\in K\cdot(a_{n-1}a_1)^*\dotsm(a_{s+1}a_{s-1})^*v.\label{b4}
\end{equation}
Finally, (\ref{b1}) and (\ref{b4}) yield
$xv=b_1\dotsm b_sb_{s+1}\dotsm b_nv\in K\cdot (a_{n-1}a_1)^*\dotsm(a_{s+1}a_{s-1})^*a_s^*v$,
which ends the proof.
\end{proof}
\end{lemma}

Now, we are ready to formulate the main result of the paper.

\begin{theorem}\label{classification}
Let $V$ be a simple left $K[C_n]$-module. Then $V$ is isomorphic to one
of the modules constructed in Proposition~\ref{simple} (in this case $n$
must be even) or $xV=0$, where $x=a_i-\lambda$ for some
$i\in\{1,\dotsc,n\}$ and $\lambda\in K$, or $x=\lambda a_j-\mu a_{j-1}$
for some $j\in\{2,\dotsc,n\}$ and $\lambda,\mu\in K$ not both equal to zero.
In the latter case $V$ may be treated as a simple left
$K[C_{n-1}]$-module and its structure can be described inductively.
\begin{proof}
Let $P$ denote the annihilator of the module $V$. Since $P$ is a
prime ideal, it follows that $P$ contains a minimal prime ideal of
$K[C_n]$, which is of the form $I_{\rho(d)}$ for some leaf $d\in
D$ (see Section~\ref{background}). We also know
that the congruence $\rho(d)$ arises as a finite extension
$\rho(d_0)\subseteq\rho(d_1)\subseteq\dotsm\subseteq\rho(d_m)=\rho(d)$,
where each $d_j$ is a diagram in level $j$ of $D$. In particular,
$I_{\rho(d_j)}\subseteq P$ for each $j=1,\dotsc,m$.

First, consider the case in which some diagram $d_{j+1}$ is
obtained from $d_j$ by adding a dot. Let us additionally assume
that $j$ is minimal with this property. If $j=0$ then the diagram $d_1$
consists of a single dot $a_i$ for some $i\in\{2,\dotsc,n-1\}$.
Hence, by (\ref{rel1-1}) and (\ref{rel1-2}) with $s=i$,
we conclude that $a_i$ is central in $K[C_n]/P$.
Therefore, by Proposition~\ref{central},
we have $xV=0$, where $x=a_i-\lambda$ for some $\lambda\in K$.
Whereas, if $j>0$ then $d_j$ consists of $j$ consecutive
arcs $\arc{a_{s+1}a_s},\dotsc,\arc{a_{s+j}a_{s-j+1}}$ and
$d_{j+1}$ arises by adjoining the dot $a_{s-j}$ to $d_j$ or by
adjoining the dot $a_{s+j+1}$ to $d_j$. If $a_{s+j}a_{s-j+1}\in P$
then Proposition~\ref{prim1} implies that $xV=0$, where
$x\in\{a_{s+j},a_{s-j+1}\}$. Whereas, if $a_{s+j}a_{s-j+1}\notin
P$ then Proposition~\ref{prim2} implies that $xV=0$, where
$x\in\{a_{s-j}-\lambda a_{s-j+1},a_{s+j+1}-\lambda
a_{s+j}:\lambda\in K\}$, and the result also follows.

Now assume that a dot does not appear in the construction of $d$,
but $d$ contains the arc $\arc{a_ja_1}$, where $j<n$ or the arc
$\arc{a_na_i}$, where $i>1$. If $a_ja_1\in P$ or $a_na_i\in P$ then
Proposition~\ref{prim1} implies that $xV=0$ for some
$x\in\{a_1,a_i,a_j,a_n\}$. Whereas, if $a_ja_1\notin P$ and
$a_na_i\notin P$ then Proposition~\ref{prim3} yields $xV=0$, where
$x\in\{a_{j+1}-\lambda a_j,a_{i-1}-\lambda a_i:\lambda\in K\}$,
hence the result follows in this situation as well.

Let us observe that if $n$ is odd then one of the cases described
above must hold. Therefore, we may assume that $n=2s$ for
some $s\ge 1$. Moreover, it remains to consider the case in which
the diagram $d$
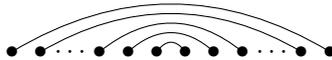
\medskip

\[\diag{
	\bullet \ar@/^1.5pc/@{-}[rrrrrrrrrrr] & \bullet \ar@/^1.2pc/@{-}[rrrrrrrrr]
	& \dotsm & \bullet \ar@/^0.9pc/@{-}[rrrrr] & \bullet \ar@/^0.6pc/@{-}[rrr]
	& \bullet \ar@/^0.3pc/@{-}[r] & \bullet & \bullet & \bullet & \dotsm & \bullet & \bullet
}\]
consists of $s$ consecutive arcs
$\arc{a_{s+1}a_s},\dotsc,\arc{a_na_1}$ (as shown in the picture).
In this situation we already know (see the proof of Lemma~\ref{lem2}) that the elements
$a_{n-j+1}a_j$ for $j=1,\dotsc,s$ are central in $K[C_n]/P$.
Therefore, by Proposition~\ref{central}, we have
$a_{n-j+1}a_j=\lambda_j$ in $K[C_n]/P$ for some $\lambda_j\in K$.
Moreover, due to Proposition~\ref{prim1}, we may assume that each
$\lambda_j\ne 0$. Further, if $a_{s+1}$ acts regularly on $V$ then
the equality $a_{s+1}(a_sa_{s+1}-\lambda_s)V=0$ implies
$(a_sa_{s+1}-\lambda_s)V=0$. Hence $a_s$ and
$\lambda_s^{-1}a_{s+1}$ are mutual inverses in $K[C_n]/P$.
Since $a_s$ commutes with $a_1,\dotsc,a_s$ in $K[C_n]/P$ (see (\ref{rel2-1})),
and $a_{s+1}$ commutes with $a_{s+1},\dotsc,a_n$ in $K[C_n]/P$ (see (\ref{rel2-2})),
we conclude that $a_s$ is a central element of $K[C_n]/P$. Thus, again
by Proposition~\ref{central}, we conclude that $xV=0$, where $x=a_s-\lambda$
for some $\lambda\in K$. Therefore, we may assume that $a_{s+1}$ does not
act regularly on $V$. In this situation Lemma~\ref{lem5}
guarantees that 
\begin{equation}\label{x0}
	\text{there exists }0\ne v\in V\text{ such that }X_0v=0
	\text{ and }a_jv=0\text{ for all }j>s
\end{equation}
(notation introduced before Lemma~\ref{lem1} is used here).
Moreover, we know that $V$ is spanned as a
$K$-linear space by elements of the set
$(a_{n-1}a_1)^*\dotsm(a_{s+1}a_{s-1})^*a_s^*v$.

First, assume $X_0=X$. We claim that in this case elements of the
set $(a_{n-1}a_1)^*\dotsm(a_{s+1}a_{s-1})^*a_s^*v$ are linearly
independent over $K$. Indeed, suppose on the contrary that
\[\sum_{i_1,\dotsc,i_s=0}^r\lambda_{i_1,\dotsc,i_s}(a_{n-1}a_1)^{i_1}
\dotsm(a_{s+1}a_{s-1})^{i_{s-1}}a_s^{i_s}v=0\]
is a nontrivial relation of linear dependence.
Then define \[m_1=\max\{i_1:\lambda_{i_1,i_2,\dotsc,i_s}\ne
0\text{ for some }i_2,\dotsc,i_s\}.\]
Observe that, by (\ref{rel2-3}) with $r=s-1$, we have
$a_1a_{n-1}a_2=a_2a_{n-1}a_1$ in $K[C_n]/P$, which yields
\[(a_na_2)(a_{n-1}a_1)=a_n(a_2a_{n-1}a_1)=a_n(a_1a_{n-1}a_2)
=(a_na_1)(a_{n-1}a_2)=\lambda_1\lambda_2\]
in $K[C_n]/P$. Since $a_na_2v=0$ and because
the element $a_na_2$ commutes in $C_n$ with $a_{n-j}a_j$ for $j>1$, we get
\begin{align*}
	0 & =(a_na_2)^{m_1}\sum_{i_1,\dotsc,i_s=0}^r\lambda_{i_1,\dotsc,i_s}(a_{n-1}a_1)^{i_1}
	\dotsm(a_{s+1}a_{s-1})^{i_{s-1}}a_s^{i_s}v\\
	&=(\lambda_1\lambda_2)^{m_1}\sum_{i_2,\dotsc,i_s=0}^r
	\lambda_{m_1,i_2,\dotsc,i_s}(a_{n-2}a_2)^{i_2}\dotsm(a_{s+1}a_{s-1})^{i_{s-1}}a_s^{i_s}v.
\end{align*}
Assume now that $k<s-1$, the numbers $m_1,\dotsc,m_k$
have already been defined, and the equality
\[\sum_{i_{k+1},\dotsc,i_s=0}^r\lambda_{m_1,\dotsc,m_k,i_{k+1},\dotsc,i_s}
(a_{n-k-1}a_{k+1})^{i_{k+1}}\dotsm(a_{s+1}a_{s-1})^{i_{s-1}}a_s^{i_s}v=0\]
holds with $\lambda_{m_1,\dotsc,m_k,i_{k+1},\dotsc,i_s}\ne 0$ for
some $i_{k+1},\dotsc,i_s$. Put
\[m_{k+1}=\max\{i_{k+1}:\lambda_{m_1,\dotsc,m_k,i_{k+1},i_{k+2},\dotsc,i_s}\ne 0
\text{ for some }i_{k+2},\dotsc,i_s\}.\]
Then, by (\ref{rel2-3}) with $r=s-k-1$, we have
$a_{k+2}a_{n-k-1}a_{k+1}=a_{k+1}a_{n-k-1}a_{k+2}$ in $K[C_n]/P$,
which yields
\begin{align*}
	(a_{n-k}a_{k+2})(a_{n-k-1}a_{k+1}) &=a_{n-k}(a_{k+2}a_{n-k-1}a_{k+1})\\
	& =a_{n-k}(a_{k+1}a_{n-k-1}a_{k+2})\\
	&=(a_{n-k}a_{k+1})(a_{n-k-1}a_{k+2})=\lambda_{k+1}\lambda_{k+2}
	\end{align*}
in $K[C_n]/P$. Since $a_{n-k}a_{k+2}v=0$ and because the element
$a_{n-k}a_{k+2}$ commutes in $C_n$ with $a_{n-j}a_j$ for $j>k+1$, we get
\begin{align*}
	0 & =(a_{n-k}a_{k+2})^{m_{k+1}}\sum_{i_{k+1},\dotsc,i_s=0}^r
	\lambda_{m_1,\dotsc,m_k,i_{k+1},\dotsc,i_s}
	(a_{n-k-1}a_{k+1})^{i_{k+1}}\dotsm(a_{s+1}a_{s-1})^{i_{s-1}}a_s^{i_s}v\\
	&=(\lambda_{k+1}\lambda_{k+2})^{m_{k+1}}\sum_{i_{k+2},\dotsc,i_s=0}^r
	\lambda_{m_1,\dotsc,m_{k+1},i_{k+2},\dotsc,i_s}
	(a_{n-k-2}a_{k+2})^{i_{k+2}}\dotsm(a_{s+1}a_{s-1})^{i_{s-1}}a_s^{i_s}v.
\end{align*}
Thus, by induction, we conclude that there exist $m_1,\dotsc,m_{s-1}$ such that
$\lambda_{m_1,\dotsc,m_{s-1},i_s}\ne 0$ for some $i_s$, and
\[\sum_{i_s=0}^r\lambda_{m_1,\dotsc,m_{s-1},i_s}a_s^{i_s}v=0.\]
Now, let $m_s=\max\{i_s:\lambda_{m_1,\dotsc,m_{s-1},i_s}\ne0\}$.
Since $a_{s+1}a_s=\lambda_s$ in $K[C_n]/P$ and $a_{s+1}v=0$, we get
\[0=a_{s+1}^{m_s}\sum_{i_s=0}^r\lambda_{m_1,\dotsc,m_{s-1},i_s}a_s^{i_s}v
=\lambda_s^{m_s}\lambda_{m_1,\dotsc,m_s}v,\]
which leads to a false conclusion that $v=0$. Therefore, the set
\[E=\{e_{i_1,\dotsc,i_s}=(\lambda_1^{-1}a_{n-1}a_1)^{i_1}\dotsm(\lambda_{s-1}^{-1}
a_{s+1}a_{s-1})^{i_{s-1}}(\lambda_s^{-1}a_s)^{i_s}v:i_1,\dotsc,i_s\ge 0\}\]
is a basis of $V$ over $K$, and one can easily check that the action of
$a_1,\dotsc,a_n\in C_n$ on the basis $E$ agrees with the action of
$a_1,\dotsc,a_n\in C_n$ on the basis of the left $K[C_n]$-module
$V(\lambda_1,\dotsc,\lambda_s)$ defined in Proposition~\ref{simple}.
Hence we get $V\cong V(\lambda_1,\dotsc,\lambda_s)$.

Finally, let us consider the last case. Namely, $X_0\ne X$. This
means that some element $a_ja_i\in X$ acts regularly on $V$ (that
is, $a_ja_i\in X_1$). Lemma~\ref{lem3} assures that we may
restrict to the situation in which $i+j=n+2$ (that is, $a_ja_i$
lies on the diagonal in the matrix notation (\ref{matrix}) of elements of $X$).
Moreover, we may assume that $j$ is minimal with that property.
In this case Lemma~\ref{lem3} and the discussion after this lemma imply that
all elements $a_la_k\in X$ with $k>i$ lie in $X_0$. Because the vector $v$
satisfies $X_0v=0$ (see (\ref{x0})), we get, in particular,
\begin{equation}\label{class-1}
	a_la_kv=0\text{ for all }a_la_k\in X\text{ with }k>i.
\end{equation}
We also know that $V$ is spanned as a $K$-linear space by the set
\[(a_{n-1}a_1)^*\dotsm(a_ja_{n-j})^*(a_{j-1}a_{n-j+1})^*(a_{j-2}a_{n-j+2})^*
\dotsm(a_{s+1}a_{s-1})^*a_s^*v.\]
Since $a_ja_i\in X_1$, Lemma~\ref{lem2} guarantees that $a_ja_i=\mu_i$
in $K[C_n]/P$ for some $0\ne\mu_i\in K$. Moreover, we have
\[(a_ja_i)(a_{j-1}a_{i-1})=a_{j-1}a_ja_ia_{i-1}=a_{j-1}a_ja_{i-1}a_i=(a_{j-1}a_i)(a_ja_{i-1})\]
in $C_n$. Since $i+j=n+2$, we get $a_{j-1}a_i=\lambda_i$ and $a_ja_{i-1}=\lambda_{i-1}$
in $K[C_n]/P$. Furthermore, we have $a_ja_i(a_{j-1}a_{i-1}a_ja_i-\lambda_{i-1}\lambda_i)V=0$.
Because $a_ja_i$ acts regularly on $V$, the last equality yields
$(a_{j-1}a_{i-1}a_ja_i-\lambda_{i-1}\lambda_i)V=0$.
Therefore, we conclude that $a_{j-1}a_{n-j+1}=a_{j-1}a_{i-1}=\mu_{i-1}$ in $K[C_n]/P$,
where $\mu_{i-1}=\lambda_{i-1}\lambda_i\mu_i^{-1}\ne 0$, hence
$V$ is also spanned as a $K$-linear space by the set
\[(a_{n-1}a_1)^*\dotsm(a_ja_{n-j})^*(a_{j-2}a_{n-j+2})^*\dotsm(a_{s+1}a_{s-1})^*a_s^*v.\]
We claim that in this case $xV=0$, where
$x=\lambda_ia_j-\mu_ia_{j-1}$. To prove this, fix
\[w=(a_{n-1}a_1)^{i_1}\dotsm(a_ja_{n-j})^{i_{n-j}}(a_{j-2}a_{n-j+2})^{i_{n-j+2}}
\dotsm(a_{s+1}a_{s-1})^{i_{s-1}}a_s^{i_s}v\in V,\]
where $i_1,\dotsc,i_{n-j},i_{n-j+2},\dotsc,i_s\ge 0$. Our aim
is to show that $xw=0$.

Assume first that $j=s+2$ (and, consequently, $i=s$). Because $a_{j-1}$ and $a_j$ commute in $C_n$
with all elements $a_{n-1}a_1,\dotsc,a_{s+2}a_{s-2}$, we have
\begin{align*}
	a_{j-1}w & =(a_{n-1}a_1)^{i_1}\dotsm(a_{s+2}a_{s-2})^{i_{s-2}}a_{j-1}a_s^{i_s}v
	\intertext{and}
	a_jw & =(a_{n-1}a_1)^{i_1}\dotsm(a_{s+2}a_{s-2})^{i_{s-2}}a_ja_s^{i_s}v.
\end{align*}
So it is enough to show that $xw'=0$, where $w'=a_s^{i_s}v$.
But $a_{j-1}v=a_jv=0$ and $a_{j-1}a_sv=\lambda_iv$, $a_ja_sv=\mu_iv$ imply that
\[a_{j-1}w'=\begin{cases}\lambda_ia_s^{i_s-1}v & \text{if }i_s>0,\\ 0 & \text{if }i_s=0\end{cases}
\qquad\text{and}\qquad a_jw'=\begin{cases}\mu_ia_s^{i_s-1}v & \text{if }i_s>0,\\ 0 & \text{if }i_s=0.\end{cases}\]
Hence the result follows in this case.

Now, let $j>s+2$.
Because $a_{j-1}$ and $a_j$ commute in $C_n$
with all elements $a_{n-1}a_1,\dotsc,a_ja_{n-j}$, we have
\begin{align*}
	a_{j-1}w & =(a_{n-1}a_1)^{i_1}\dotsm(a_ja_{n-j})^{i_{n-j}}a_{j-1}
	(a_{j-2}a_{n-j+2})^{i_{n-j+2}}\dotsm(a_{s+1}a_{s-1})^{i_{s-1}}a_s^{i_s}v
	\intertext{and}
	a_jw & =(a_{n-1}a_1)^{i_1}\dotsm(a_ja_{n-j})^{i_{n-j}}a_j
	(a_{j-2}a_{n-j+2})^{i_{n-j+2}}\dotsm(a_{s+1}a_{s-1})^{i_{s-1}}a_s^{i_s}v.
\end{align*}
So it suffices to check that $xw'=0$, where
$w'=(a_{j-2}a_{n-j+2})^{i_{n-j+2}}\dotsm(a_{s+1}a_{s-1})^{i_{s-1}}a_s^{i_s}v$.
Suppose that $i_{n-j+2}>0$. Then, remembering that $i+j=n+2$, we get
\begin{align*}
	a_{j-1}w' & =a_{j-1}a_{j-2}a_{n-j+2}(a_{j-2}a_{n-j+2})^{i_{n-j+2}-1}
	(a_{j-3}a_{n-j+3})^{i_{n-j+3}}\dotsm(a_{s+1}a_{s-1})^{i_{s-1}}a_s^{i_s}v\\
	& =(a_{j-1}a_{n-j+2})a_{j-2}(a_{j-2}a_{n-j+2})^{i_{n-j+2}-1}
	(a_{j-3}a_{n-j+3})^{i_{n-j+3}}\dotsm(a_{s+1}a_{s-1})^{i_{s-1}}a_s^{i_s}v\\
	& =\lambda_ia_{j-2}(a_{j-2}a_{n-j+2})^{i_{n-j+2}-1}
	(a_{j-3}a_{n-j+3})^{i_{n-j+3}}\dotsm(a_{s+1}a_{s-1})^{i_{s-1}}a_s^{i_s}v\\
	\intertext{and}
	a_jw' & =a_ja_{j-2}a_{n-j+2}(a_{j-2}a_{n-j+2})^{i_{n-j+2}-1}
	(a_{j-3}a_{n-j+3})^{i_{n-j+3}}\dotsm(a_{s+1}a_{s-1})^{i_{s-1}}a_s^{i_s}v\\
	& =(a_ja_{n-j+2})a_{j-2}(a_{j-2}a_{n-j+2})^{i_{n-j+2}-1}
	(a_{j-3}a_{n-j+3})^{i_{n-j+3}}\dotsm(a_{s+1}a_{s-1})^{i_{s-1}}a_s^{i_s}v\\
	& =\mu_ia_{j-2}(a_{j-2}a_{n-j+2})^{i_{n-j+2}-1}
	(a_{j-3}a_{n-j+3})^{i_{n-j+3}}\dotsm(a_{s+1}a_{s-1})^{i_{s-1}}a_s^{i_s}v,
\end{align*}
because $a_{j-1}a_{n-j+2}=a_{j-1}a_i=\lambda_i$ and
$a_ja_{n-j+2}=a_ja_i=\mu _i$ in $K[C_n]/P$.
Hence $xw'=0$ in this case. Next, assume that $i_{n-j+2}=0$.
If all $i_{n-j+3}=\dotsb=i_{s-1}=0$ then $w'=a_s^{i_s}v$.
Since $j>s+2$, we have $a_{j-1}a_s,a_ja_s\in X$. 
Thus (\ref{class-1}) gives $a_{j-1}a_sv=a_ja_sv=0$, because $s>i$.
Moreover, we have $a_{j-1}v=a_jv=0$. Therefore, $a_{j-1}w'=a_jw'=0$
and, in consequence, $xw'=0$. Finally, assume that
$i_{n-j+2}=0$ but $i_k>0$ for some $k\in\{n-j+3,\dotsc,s-1\}$,
and choose minimal $k$ with this property. In this situation we have
$w'=(a_{n-k}a_k)^{i_k}\dotsm(a_{s+1}a_{s-1})^{i_{s-1}}a_s^{i_s}v$.
Because $k\ge n-j+3$ then both $j-1$ and $j$ are $\ge n-k$, hence
$a_{j-1}a_{n-k}a_k=a_{n-k}(a_{j-1}a_k)$ and $a_ja_{n-k}a_k=a_{n-k}(a_ja_k)$ in $C_n$.
Therefore
\begin{align*}
	a_{j-1}w' & =a_{j-1}a_{n-k}a_k(a_{n-k}a_k)^{i_k-1}(a_{n-k-1}a_{k+1})^{i_{k+1}}
	\dotsm(a_{s+1}a_{s-1})^{i_{s-1}}a_s^{i_s}v,\\
	& =a_{n-k}(a_{j-1}a_k)(a_{n-k}a_k)^{i_k-1}(a_{n-k-1}a_{k+1})^{i_{k+1}}
	\dotsm(a_{s+1}a_{s-1})^{i_{s-1}}a_s^{i_s}v,\\	
	& =a_{n-k}(a_{n-k}a_k)^{i_k-1}a_{j-1}a_k(a_{n-k-1}a_{k+1})^{i_{k+1}}
	\dotsm(a_{s+1}a_{s-1})^{i_{s-1}}a_s^{i_s}v
	\intertext{and}	
	a_jw' & =a_ja_{n-k}a_k(a_{n-k}a_k)^{i_k-1}(a_{n-k-1}a_{k+1})^{i_{k+1}}
	\dotsm(a_{s+1}a_{s-1})^{i_{s-1}}a_s^{i_s}v.\\
	& =a_{n-k}(a_ja_k)(a_{n-k}a_k)^{i_k-1}(a_{n-k-1}a_{k+1})^{i_{k+1}}
	\dotsm(a_{s+1}a_{s-1})^{i_{s-1}}a_s^{i_s}v.\\
	& =a_{n-k}(a_{n-k}a_k)^{i_k-1}a_ja_k(a_{n-k-1}a_{k+1})^{i_{k+1}}
	\dotsm(a_{s+1}a_{s-1})^{i_{s-1}}a_s^{i_s}v.
\end{align*}
Further, observe that $k\ge n-j+3$ implies that for each $k<l<s$ we have $j-1>n-l$, hence
$a_{j-1}a_k$ and $a_ja_k$ commute with $a_{n-l}a_l$ in $C_n$.
Since $k<s<j-1$ it is also clear that $a_{j-1}a_k$ and $a_ja_k$ commute in $C_n$ with $a_s$.
Thus
\begin{align*}
	a_{j-1}w' & =a_{n-k}(a_{n-k}a_k)^{i_k-1}(a_{n-k-1}a_{k+1})^{i_{k+1}}
	\dotsm(a_{s+1}a_{s-1})^{i_{s-1}}a_s^{i_s}a_{j-1}a_kv,\\
	a_jw' & =a_{n-k}(a_{n-k}a_k)^{i_k-1}(a_{n-k-1}a_{k+1})^{i_{k+1}}
	\dotsm(a_{s+1}a_{s-1})^{i_{s-1}}a_s^{i_s}a_ja_kv.
\end{align*}
Since $a_{j-1}a_kv=a_ja_kv=0$ (by (\ref{class-1}), because $k\ge n-j+3=i+1$ implies
that $a_{j-1}a_k,a_ja_k\in X$ and $k>i$),
we get $a_{j-1}w'=a_jw'=0$ and, in consequence, $xw'=0$. This finishes the proof.
\end{proof}
\end{theorem}

Recall that a representation of a monoid $M$ in a $K$-linear space
$V$ is said to be monomial, if $V$ admits a basis $E$ such that
for each $w\in M$ and each $e\in E$ there exist $\lambda\in K$ and
$f\in E$ such that $we=\lambda f$. As a consequence of
Proposition~\ref{simple} and Theorem~\ref{classification} we get
the following remarkable result.

\begin{corollary}
Each irreducible representation of the Chinese monoid $C_n$ is monomial.
\end{corollary}

This is in contrast with the results obtained for the, similarly defined,
important class of plactic algebras. Namely, in \cite{cedo-kubat-okn}
it is shown that the plactic algebra of rank $4$ admits irreducible
representations which are not monomial. It is also worth to note that
all irreducible representations of plactic algebras of rank
not exceeding $3$ are monomial (see \cite{ko1}).

\section{\texorpdfstring{Illustration of the main theorem for $n\le 4$}
{Illustration of the main theorem for $n<5$}}\label{illustration}
In order to provide more insight into the nature of
Theorem~\ref{classification}, we interpret it in the case of small
values of $n$. The case $n=1$ is trivial. Next, it is well known
that the Chinese algebra $K[C_2]$ of rank $2$ coincides with the
plactic algebra of rank $2$. Moreover, the irreducible
representations of $C_2$ are easy to describe, as they are induced
from irreducible representations of the bicyclic monoid
$B\cong C_2/(a_2a_1=1)$. Namely, we have the following result
(see \cite{ko1} for more details).
\begin{remark}\label{rem1}
    Let $V$ be a simple left $K[C_2]$-module. Then $V$ is $1$-dimensional
    or $V\cong Z$, where $Z$ is the simple left $K[C_2]$-module defined
    just before Proposition~\ref{inductive}.
\end{remark}
Our next step is to describe all irreducible representations of the monoid $C_3$.
In this case the diagram $D$ has the form\vspace{-12pt}

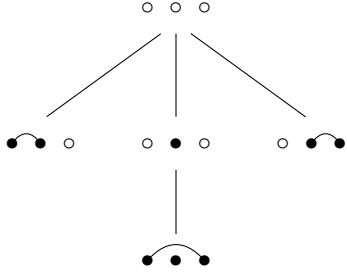
\newbox{\zero}
\newbox{\oneD}
\newbox{\oneAL}
\newbox{\oneAR}
\newbox{\two}

\savebox{\zero}{\begin{varwidth}[t][-0.5cm][c]{\textwidth}
$\diag{ \circ & \circ & \circ }$\end{varwidth}}
\savebox{\oneD}{\begin{varwidth}[t][0.5cm][c]{\textwidth}
$\diag{ \circ & \bullet & \circ }$\end{varwidth}}
\savebox{\oneAL}{\begin{varwidth}[t][0.5cm][c]{\textwidth}
$\diag{ \bullet \ar@/^0.3pc/@{-}[r] & \bullet & \circ }$\end{varwidth}}
\savebox{\oneAR}{\begin{varwidth}[t][0.5cm][c]{\textwidth}
$\diag{ \circ & \bullet \ar@/^0.3pc/@{-}[r]& \bullet }$\end{varwidth}}
\savebox{\two}{\begin{varwidth}[t][0.5cm][c]{\textwidth}
$\diag{ \bullet \ar@/^0.5pc/@{-}[rr] & \bullet & \bullet }$\end{varwidth}}

\[\xymatrix{
& \usebox{\zero}  \ar@{-}[ld]  \ar@{-}[d] \ar@{-}[rd] \\
\usebox{\oneAL} & \usebox{\oneD} \ar@{-}[d] & \usebox{\oneAR} \\
& \usebox{\two}}\]
and three leaves of this diagram correspond to the minimal prime ideals
of $K[C_3]$:
\begin{align*}
	P_1 &= (a_2,a_3\text{ commute},\,a_2a_1\text{ central}),\\
	P_2 & =(a_2\text{ central}),\\
	P_3 & = (a_1,a_2\text{ commute},\,a_3a_2\text{ central}).
\end{align*}
Here, writing for example `$a_2,a_3$ commute'
in $P_1$ we mean that $P_1$ contains the element $a_2a_3-a_3a_2$.
Similarly, writing `$a_2a_1$ central' we understand that $P_1$ contains
all elements of the form $a_ia_2a_1-a_2a_1a_i$ for $i=1,2,3$.
The same convention applies to other minimal prime ideals of $K[C_3]$.

Hence, by Remark~\ref{rem1} and the results from
Section~\ref{irr}, we get the following classification.
\begin{remark}\label{rem2}
Let $V$ be a simple left $K[C_3]$-module. Then $V$ is $1$-dimensional or there
exists a basis $\{e_i:i\ge 0\}$ of $V$ such that exactly one of the
following possibilities holds:
\begin{enumerate}
	\item there exist $\lambda,\mu\in K$ such that $\lambda\ne 0$ and
	$a_1e_i=\lambda e_{i+1}$, $a_2e_i=e_{i-1}$, $a_3e_i=\mu e_{i-1}$ for all $i\ge 0$.
	\item there exist $\lambda,\mu\in K$ such that $\lambda\ne 0$ and
	$a_1e_i=\lambda e_{i+1}$, $a_2e_i=\mu e_i$, $a_3e_i=e_{i-1}$ for all $i\ge 0$.
	\item there exist $\lambda,\mu\in K$ such that $\mu\ne 0$ and
	$a_1e_i=\lambda e_{i+1}$, $a_2e_i=\mu e_{i+1}$, $a_3e_i=e_{i-1}$ for all $i\ge 0$.
\end{enumerate}
Note that, to make our statements more compact, we adopted the
convention that $e_{-1}=0$.
\end{remark}
Finally, let us describe all irreducible representations of the monoid $C_4$.
In this situation the diagram $D$ has the form\vspace{-12pt}

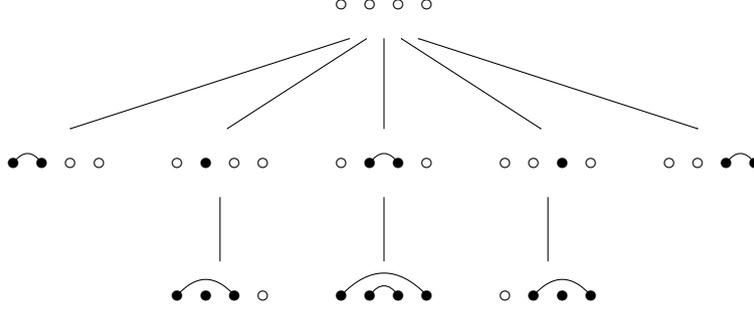
\newbox{\zero}
\newbox{\oneDL}
\newbox{\oneDR}
\newbox{\oneAL}
\newbox{\oneAC}
\newbox{\oneAR}
\newbox{\twoADL}
\newbox{\twoADR}
\newbox{\twoAA}

\savebox{\zero}{\begin{varwidth}[t][-0.7cm][c]{\textwidth}
$\diag{ \circ & \circ & \circ & \circ}$\end{varwidth}}
\savebox{\oneDL}{\begin{varwidth}[t][0.7cm][c]{\textwidth}
$\diag{ \circ & \bullet & \circ & \circ}$\end{varwidth}}
\savebox{\oneDR}{\begin{varwidth}[t][0.7cm][c]{\textwidth}
$\diag{ \circ & \circ & \bullet & \circ}$\end{varwidth}}
\savebox{\oneAL}{\begin{varwidth}[t][0.7cm][c]{\textwidth}
$\diag{ \bullet \ar@/^0.3pc/@{-}[r] & \bullet & \circ & \circ}$\end{varwidth}}
\savebox{\oneAC}{\begin{varwidth}[t][0.7cm][c]{\textwidth}
$\diag{ \circ & \bullet \ar@/^0.3pc/@{-}[r]& \bullet & \circ}$\end{varwidth}}
\savebox{\oneAR}{\begin{varwidth}[t][0.7cm][c]{\textwidth}
$\diag{ \circ & \circ & \bullet \ar@/^0.3pc/@{-}[r] & \bullet}$\end{varwidth}}
\savebox{\twoADL}{\begin{varwidth}[t][0.7cm][c]{\textwidth}
$\diag{ \bullet \ar@/^0.5pc/@{-}[rr] & \bullet & \bullet & \circ}$\end{varwidth}}
\savebox{\twoADR}{\begin{varwidth}[t][0.7cm][c]{\textwidth}
$\diag{ \circ & \bullet \ar@/^0.5pc/@{-}[rr] & \bullet & \bullet } $ \end{varwidth}}
\savebox{\twoAA}{\begin{varwidth}[t][0.7cm][c]{\textwidth}
$\diag{ \bullet \ar@/^0.7pc/@{-}[rrr] & \bullet \ar@/^0.3pc/@{-}[r] & \bullet & \bullet } $\end{varwidth}}
 
\[\xymatrix{
	&& \usebox{\zero} \ar@{-}[lld] \ar@{-}[ld]  \ar@{-}[d] \ar@{-}[rd] \ar@{-}[rrd] \\
	\usebox{\oneAL}  & \usebox{\oneDL} \ar@{-}[d] & \usebox{\oneAC} \ar@{-}[d]
	& \usebox{\oneDR} \ar@{-}[d] & \usebox{\oneAR} \\
	& \usebox{\twoADL} & \usebox{\twoAA} & \usebox{\twoADR}
}\]
and five leaves of this diagram correspond to the minimal prime ideals of $K[C_4]$:
\begin{align*}
	P_1 & = (a_2,a_3,a_4\text{ commute},\,a_2a_1,a_3a_1\text{ central}),\\
	P_2 & = (a_3,a_4\text{ commute},\,a_2,a_3a_1\text{ central}),\\
	P_3 & = (a_1,a_2\text{ commute},\,a_3,a_4\text{ commute},\,a_3a_2\text{ central}),\\
	P_4 & = (a_1,a_2\text{ commute},\,a_3,a_4a_2\text{ central}),\\
	P_5 & = (a_1,a_2,a_3\text{ commute},\,a_4a_2,a_4a_3\text{ central}).
\end{align*}

Now, Remark~\ref{rem2} together with the results
obtained in Section~\ref{irr} lead to the following
classification.
\begin{remark}
Let $V$ be a simple left $K[C_4]$-module. Then $V$ is $1$-dimensional or
there exists a basis $\{e_{i,j}:i,j\ge 0\}$ of $V$ and $0\ne\lambda,\mu\in K$ such that
\[a_1e_{i,j}=\lambda e_{i+1,j+1},\qquad a_2e_{i,j}=\mu e_{i,j+1},\qquad
a_3e_{i,j}=e_{i,j-1},\qquad a_4e_{i,j}=e_{i-1,j-1}\] for all $i,j\ge 0$
(with the convention that $e_{i,j}=0$ if $i=-1$ or $j=-1$), or there exists a basis
$\{e_i:i\ge 0\}$ of $V$ such that exactly one of the following possibilities holds:
\begin{enumerate}
	\item[(1)] there exist $\lambda,\mu,\nu\in K$ such that $\lambda\ne 0$ and
	$a_1e_i=\lambda e_{i+1}$, $a_2e_i=e_{i-1}$, $a_3e_i=\mu e_{i-1}$, $a_4e_i=\nu e_{i-1}$ for all $i\ge 0$.
	\item[(2)] there exist $\lambda,\mu,\nu\in K$ such that $\lambda\ne 0$ and
	$a_1e_i=\lambda e_{i+1}$, $a_2e_i=\mu e_i$, $a_3e_i=e_{i-1}$, $a_4e_i=\nu e_{i-1}$ for all $i\ge 0$.
	\item[(3.1)] there exist $\lambda,\mu,\nu\in K$ such that $\mu\ne 0$ and
	$a_1e_i=\lambda e_{i+1}$, $a_2e_i=\mu e_{i+1}$, $a_3e_i=e_{i-1}$, $a_4e_i=\nu e_{i-1}$ for all $i\ge 0$.
	\item[(3.2)] there exist $\lambda,\mu,\nu\in K$ such that $\lambda\ne 0$ but $\mu\nu=0$ and
	$a_1e_i=\lambda e_{i+1}$, $a_2e_i=\mu e_i$, $a_3e_i=\nu e_i$, $a_4e_i=e_{i-1}$ for all $i\ge 0$.
	\item[(4)] there exist $\lambda,\mu,\nu\in K$ such that $\mu\ne 0$ and
	$a_1e_i=\lambda e_{i+1}$, $a_2e_i=\mu e_{i+1}$, $a_3e_i=\nu e_i$, $a_4e_i=e_{i-1}$ for all $i\ge 0$.
	\item[(5)] there exist $\lambda,\mu,\nu\in K$ such that $\nu\ne 0$ and
	$a_1e_i=\lambda e_{i+1}$, $a_2e_i=\mu e_{i+1}$, $a_3e_i=\nu e_{i+1}$, $a_4e_i=e_{i-1}$ for all $i\ge 0$.
\end{enumerate}
Note that, as in Remark~\ref{rem2}, we used the convention that $e_{-1}=0$. Moreover, it is worth to notice
that modules in family $(i)$, for $i=1,2,4,5$, contain in their annihilators the ideal $P_i$.
Furthermore, modules in both families (3.1) and (3.2) contain $P_3$ in their annihilators.
\end{remark}

\end{document}